\documentclass[VANCOUVER]{WileyNJD-v1}
\articletype{Research Article}%

\usepackage{amsfonts}
\usepackage{amssymb}
\usepackage{amsmath}
\usepackage{graphicx}
\usepackage{epstopdf}
\usepackage{float}
\usepackage{bbold}
\usepackage{tikz}
\usepackage{amsopn}
\usepackage{hyperref}
\usepackage{cleveref}

\usetikzlibrary{arrows,matrix,positioning}
\ifpdf
\DeclareGraphicsExtensions{.eps,.pdf,.png,.jpg}
\else
\DeclareGraphicsExtensions{.eps}
\fi

\newcommand{\Sopt}{S_{\mathrm{opt}}}
\newcommand{\Sreal}{S_{\mathrm{real}}}
\usepackage{ulem} 
\usepackage{color}

\begin{document}

\title{Computing the minimal rebinding effect for non-reversible processes\protect}

\author[1]{Susanne R\"ohl*}

\author[1]{Marcus Weber}

\author[2,1]{Konstantin Fackeldey}

\authormark{Susanne R\"ohl \textsc{et al}}

\address[1]{\orgname{Zuse Institute Berlin (ZIB)}, \orgaddress{\state{Takustra{\ss}e 7, 14195 Berlin}, \country{Germany}}}

\address[2]{\orgdiv{Institut f\"ur Mathematik}, \orgname{Technische Universit\"at Berlin}, \orgaddress{\state{Stra{\ss}e des 17. Juni 136, 10623 Berlin}, \country{Germany}}}

\corres{*Susanne R\"ohl, Zuse Institute Berlin, Takustra{\ss}e 7, 14195 Berlin. \email{susanne.roehl@fu-berlin.de}}


\abstract[Abstract]{The aim of this paper is to investigate the rebinding effect, a phenomenon describing a ``short-time memory'' which can occur when projecting a Markov process onto a smaller state space.
	For guaranteeing a correct mapping by the Markov State Model, we assume a fuzzy clustering in terms of membership functions, assigning degrees of membership to each state.
	The macro states are represented by the membership functions and may be overlapping.
	The magnitude of this overlap is a measure for the strength of the rebinding effect, caused by the projection and stabilizing the system.
	A minimal bound for the rebinding effect included in a given system is computed as the solution of an optimization problem.
	Based on membership functions chosen as a linear combination of Schur vectors we are thus able to
	to compute the minimal rebinding effect for non-reversible processes.}



\maketitle


\section{Introduction}\label{sec1}

Markov processes are memoryless stochastic processes with applications in many different kinds of areas.
They are employed to describe molecular systems like protein folding\cite{da2014application} or ligand-binding processes\cite{shan2011}.
Such processes act on high dimensional state spaces and additionally require simulations on rather long time-scales in order to observe rare conformational changes. 
Consequently, a reduction of dimension is aimed at, which can be realized by a projection onto a smaller state space.
The reduced model should represent the correct long-time behaviour of the process, while being less complex.
The existence of metastable sets can be exploited to create such a ``Markov State Model''\cite{bowman2013introduction, chodera2014markov,milestoning}.

A well-established solution is the  fuzzy clustering algorithm PCCA+, which identifies \nopagebreak metastable sets with the aid of membership functions $\chi = XA$, being a linear combination of eigenvectors\cite{weber2006meshless}.

When projecting a process onto a finite state space, it may lose its Markov property, more precisely it can include short-time memory effects.
Such memory effects were detected in the context of ligand-binding-systems, where in certain configurations significantly increased binding affinities were observed\cite{vauquelin2010}.
They are explained by an additional memory caused by the projection: short time after a ligand unbounds from its target, it is assumed to be still nearby and thus rebinds with a high probability. Consequently, this short-time memory is denoted as \textbf{rebinding effect}.
This memory effect is strongly related to the \textbf{overlap} of the membership functions $\chi$ determining the clustering.
Hence, knowing them makes it easy to compute the actual rebinding effect caused by this projection.
However, in many cases the original process and the membership functions are not known. For instance, a finite process can be constructed as the solution of a differential equation and just be interpreted as the projection of a larger process. In order to identify possible memory effects included in that system, it is favorable to estimate the rebinding effect. This can be achieved by solving an optimization problem, revealing a minimal bound:
Given a clustered system, how much rebinding is included \textbf{at least}?

The computation of the minimal rebinding effect included in a given kinetics has been accomplished for reversible processes in 2014 by Weber and Fackeldey\cite{weber2014}.
In this paper, the formulation of the corresponding optimization problem is extended onto non-reversible processes.
This is achieved by employing the framework of GenPCCA, a recent modification of PCCA+ by Fackeldey and Weber\cite{fackeldey2017gen} from 2017, which is based on Schur vectors instead of eigenvectors and includes non-reversible processes.
This generalization is of particular interest since many real-world processes are non-reversible\cite{fackeldey2017}.
\\

A significant application of the presented topic lies in the area of computational drug design.
In order to treat diseases, ligands are designed such that they bind to pathogenic target molecules.
Improving the binding affinity is one important goal in drug design.
For a precise prediction of the binding affinity, possible rebinding effects need to be considered, since they can influence the binding behaviour.
\\

The article is organized as follows.
In \cref{sec:projection}, we introduce the physical and mathematical framework which is necessary to describe the time-evolution of molecular systems and their projections onto finite spaces.
For that purpose, the concept of a real Schur decomposition plays an important role and different possible shapes of such a decomposition will be analyzed.
Afterwards, we describe the rebinding effect in the context of a receptor-ligand system and set in relation to the choice of the projection.
In \cref{sec:optimization}, we present an optimization problem providing a lower bound for the rebinding effect included in a given molecular kinetics, which is valid for reversible \textbf{and} nonreversible processes.
Finally, we validate the results on some illustrative examples in \cref{sec:numerical}.

\section{Projection of a Molecular System}
\label{sec:projection}

A molecular system consisting of $N$ particles can be represented in a $6N$-dimensional phase space $\Gamma = \Omega \times \mathbb{R}^{3N}$, including the position and momentum coordinates of all particles.
Since conformational changes are of particular interest, such a system is usually described by a continuous transfer operator acting on configuration space $\Omega$, see e.g. \cite{schutte2001transfer,weber2011subspace}.
However, instead of considering the continuous case, we start directly with a discretized version acting on an $m$-dimensional state space $E = \{1, \dots, m\}$, being a subset of the configuration space. This process is characterized by a finite transition matrix $P := P(\tau) \in \mathbb{R}^{m \times m}$ and a stationary distribution $\pi \in \mathbb{R}^m$, which is assumed to be unique.

The micro states will be clustered conveniently, such that the resulting macro states represent the metastable conformations of the molecular system. For considering non-reversible processes, the Schur decomposition is of particular importance.
In the following, we briefly summarize the mathematical concepts for these two main topics.

\subsection{Fuzzy Clustering}

Let $1 = \lambda_1 > | \lambda_2 | \geq \dots \geq | \lambda_n|$ be the \textbf{dominant} spectrum of the transition matrix $P$, i.e. the eigenvalues of largest absolute value which are well-separated from the rest of the spectrum.
Let $X = \{ X_1, \dots, X_n\}$ be a matrix of associated real orthogonal Schur vectors, i.e. vectors fulfilling $PX = X \Lambda$, where $\Lambda$ is a \textbf{real} Schur decomposition. Then $\Lambda$ is of \textbf{block}-triagonal shape and has $\lambda_1, \dots, \lambda_n$ as eigenvalues.
According to GenPCCA\cite{fackeldey2017gen}, membership functions $\chi_1, \dots, \chi_n: E \rightarrow [0,1]$ can be built as a linear combination
\begin{equation*}
	\chi = XA
\end{equation*}
of the dominant Schur vectors with a regular transformation matrix $A \in \mathbb{R}^{n \times n}$.
Let $\langle \cdot,\cdot \rangle_\pi$ be the $\pi$-weighted $L^2$ scalar product, by using a Galerkin projection, this choice of membership functions yields a matrix representation
\begin{equation}
\label{eq:matrix_representation}
P_c (\tau) = S^{-1} T = \langle \chi, \chi \rangle_\pi^{-1} \langle \chi, P(\tau) \chi \rangle_\pi
\end{equation}
with two stochastic matrices $S$ and $T$. They are given by
\begin{equation}
\label{eq:matrix_representation2}
\begin{aligned}
T & = D^{-1} \langle \chi, P (\tau) \chi \rangle_\pi = D^{-1}A^T \Lambda A \ \ \textrm{ and } \\
S & = D^{-1} \langle \chi, \chi \rangle_\pi = D^{-1} A^T A,
\end{aligned}
\end{equation}
where $D = \mathrm{diag}(w_1, \dots, w_n)$ is the diagonal matrix consisting of the statistical weights
\begin{equation*}
	w_i = \langle \chi_i, \mathbb 1 \rangle_\pi
\end{equation*}
of the conformations $\chi_i$, i.e. the probabilities of the clustered process to be in the conformations in equilibrium.
The clustered process can according to \eqref{eq:matrix_representation} and \eqref{eq:matrix_representation2} as well be represented by
\begin{equation}
\label{eq:matrix_representation3}
P_c = A^{-1} \Lambda A.
\end{equation}
The transition rate matrix $Q$ is related to the transition matrix via $\mathrm{exp}(\tau Q) = P(\tau)$.
The projection of a discretized transition rate matrix $Q \in \mathbb{R}^{m \times m}$ works similar to the above procedure and yields a matrix representation
\begin{equation}
\label{eq:matrix_representation4}
Q_c = A^{-1} \Xi A = \langle \chi, \chi \rangle_\pi^{-1} \langle \chi, Q\chi \rangle_\pi
\end{equation}
with the real Schur decomposition $\Xi$ corresponding to the $n$ dominant eigenvalues $0 = \xi_1, \xi_2 \dots, \xi_n$ of $Q$ and $A \in \mathbb{R}^{n \times n}$ the transformation matrix obtained by GenPCCA\cite{fackeldey2017gen}, providing an \textbf{optimal} solution.
The eigenvalues of the transition matrix and the transition rate matrix are related via
\begin{equation}
\mathrm{exp}(\xi_i) = \lambda_i.
\end{equation}
The Schur decomposition of a reversible process is equal to its spectral decomposition.
In that case, the Schur matrices $\Lambda, \Xi$ are \textbf{diagonal} matrices consisting of the real eigenvalues $1 = \lambda_1 > \dots \geq \lambda_n$ of $P$ and $0 = \xi_1 > \dots \geq \xi_n$ of $Q$.
In contrast to the well-known clustering algorithm PCCA+\cite{deuflhard2005robust} providing a solution only for reversible processes, the generalized version GenPCCA includes reversible as well as non-reversible processes.
Apart from the fact that GenPCCA takes Schur vectors instead of eigenvectors as input, the algorithm remains the same.

\subsection{Reversibility -- Non-reversibility}

A Markov chain given by the transition matrix $P \in \mathbb{R}^{m \times m}$ is reversible, if detailed balance is fulfilled, i.e. if the matrix $DP$ is symmetric. In this case, the diagonal matrix $D = \mathrm{diag}(\pi_1, \dots, \pi_m)$ consists of the entries of the stationary distribution $\pi =(\pi_1 \dots\pi_m)^T \in \mathbb{R}^m$, i.e. $\pi$ is a vector which meets $\pi^T P = \pi^T$.

In contrast to most of the existing literature, we employ a real Schur decomposition instead of the spectral decomposition for the clustering, because we are interested in an investigation including reversible as well as non-reversible processes.
In order to be feasible, the above presented algorithm requires real and orthogonal vectors spanning an invariant subspace.
Even though this procedure works for reversible processes using a set of dominant eigenvectors, the requirements are not necessarily fulfilled for non-reversible processes. Some problems that can occur:
\begin{itemize}
	\item $P$ has real eigenvalues, but non-orthogonal eigenvectors,
	\item $P$ is non-diagonalizable,
	\item $P$ has complex eigenvalues, leading to complex eigenvectors.
\end{itemize}

Since reversibility of a process cannot be presumed for real-world processes (e.g. measuring errors), we employ the generalized approach in terms of a real Schur decomposition instead of the spectral decomposition. This approach avoids the aforementioned problems: the real Schur decomposition exists for all transition matrices $P$ and yields a set of real and orthogonal Schur vectors.
We will be able to exploit the structure of the real Schur decomposition
\begin{equation*}
	\Lambda = 
	\begin{pmatrix}
		A_1 & * & \cdots & * \\
		0 & A_2 & \cdots & * \\
		\vdots & \vdots & \ddots & \vdots \\
		0 & 0 & \cdots & A_n
	\end{pmatrix}.
\end{equation*}
induced by the spectrum of $P$. We obtain a block-triagonal shape with blocks $A_i$ s.t. each $1 \times 1$-block corresponds to a real eigenvalue and each $2 \times 2$-block corresponds to a pair of complex conjugate eigenvalues.

Even though considering a decomposition with ordered blocks, e.g. according to \cite{brandts2002matlab}, we have to bear in mind that the real Schur decomposition is \textbf{not} unique. This has been discussed in \cite{ReuterFackeldeyWeber}.

\section{Rebinding Effect}
\label{sec:rebinding}

The projection of a molecular system on a smaller state space can lead to a short-time memory included in the clustered process.
This phenomena can occur in all kinds of processes when projecting them. We introduce it on an easy example and show how it can be measured employing the mathematical tools from section \ref{sec:projection}.

\subsection{Mathematical Model of Receptor-Ligand System}

The binding behaviour of a simple receptor-ligand system is formalized as follows.
A ligand (L) can bind to a receptor (R) and form a receptor-ligand complex (LR) which can dissociate again into its original components. This process can be represented by a reaction equation
\begin{equation}
\label{eq:reaction}
\mathrm{L} + \mathrm{R} 	\overset{k_{\mathrm{on}}}{\underset{k_{\mathrm{off}}}{\rightleftharpoons}}	 \mathrm{LR}.
\end{equation}
Being a process in chemical equilibrium, the law of mass action states that the ratio between the concentration of reactants and products is constant. 
The corresponding \textit{dissociation constant} $k_d$ is given by
\begin{equation*}
	k_d
	= \frac{k_{\mathrm{off}}}{k_{\mathrm{on}}}
	= \frac{\mathrm{[L]} \cdot \mathrm{[R]}}{\mathrm{[LR]}},
\end{equation*}
where [L] represents the concentration of unbound ligands, [R] the concentration of unoccupied receptors and [LR] the concentration of receptor-ligand complexes, respectively.
This constant is used to describe the \textit{binding affinity} between a ligand and a receptor, that is how strongly the ligand can bind to his particular receptor. If the dissociation constant is small, then there are relatively many complexes in comparison to unbound molecules, and for this reason, the binding affinity between the ligand and the receptor is high.
The \textit{association constant} $k_a$ is the inverse of the dissociation constant
\begin{equation*}
	k_a
	= \frac{k_{\mathrm{on}}}{k_{\mathrm{off}}}
	= \frac{\mathrm{[LR]}}{\mathrm{[L]} \cdot \mathrm{[R]}}.
\end{equation*}
There are different factors which can influence the binding affinity of a process.
It depends on the nature of the constituent molecules, like their shape, size and possible charge.
The binding affinity of a particular ligand-protein interaction can also significantly change with solution conditions, e.g. temperature, pH or salt concentration.
For instance, a higher temperature leads to a faster movement of the molecules and therefore increases the probability of binding events.
In general, high-affinity binding results in a higher degree of occupancy of the receptors than it is the case for low-affinity binding; the residence time does not correlate\cite{lauffenburger1993receptors}.
\\

Starting from the reaction equation \eqref{eq:reaction}, we claim that a ligand can be found in two different macro states: ``unbound'' (L) or ``bound'' (LR). Then the probabilities of the ligand to be in one of these states are described by the probability vector $x^T = \frac{1}{s}(\mathrm{[L],[LR]})$, where $s = \mathrm{[L]} + \mathrm{[LR]} = \textrm{const.}$ is the normalization constant.
This leads to an ordinary differential equation
\begin{equation*}
	\dot{x}^T = x^T Q_c.
\end{equation*}
The matrix $Q_c$ consists of the rates of reaction,
\begin{equation}
\label{eq:model_monovalent}
Q_c = 
\begin{pmatrix}
-k_a[R] & k_a[R]  \\
k_d      & -k_d
\end{pmatrix},
\end{equation}
where $k_a$ and $k_d$ are the association and dissociation constants. It corresponds to the transition rate matrix of a Markov chain, that means it describes a \textbf{memoryless} process. \\

\begin{minipage}{\linewidth}
	\centering
	\begin{minipage}[t]{0.2\textwidth}
		\includegraphics[width=\textwidth]{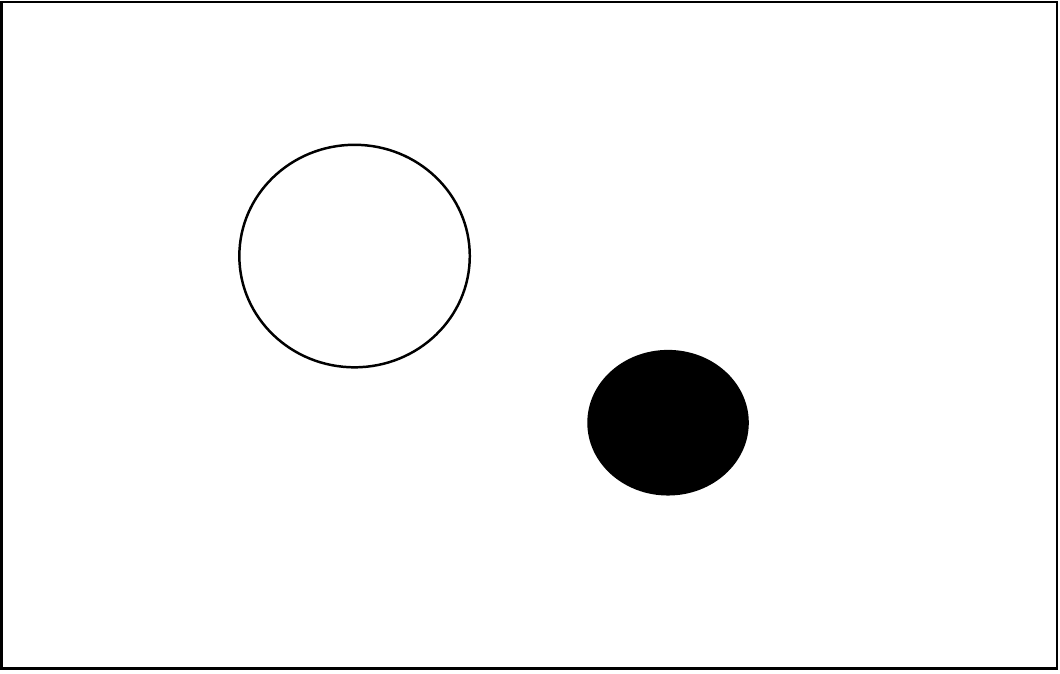}
		\centering {\small (a) ``unbound''}
	\end{minipage}
	\begin{minipage}[t]{0.2\textwidth}
		\includegraphics[width=\textwidth]{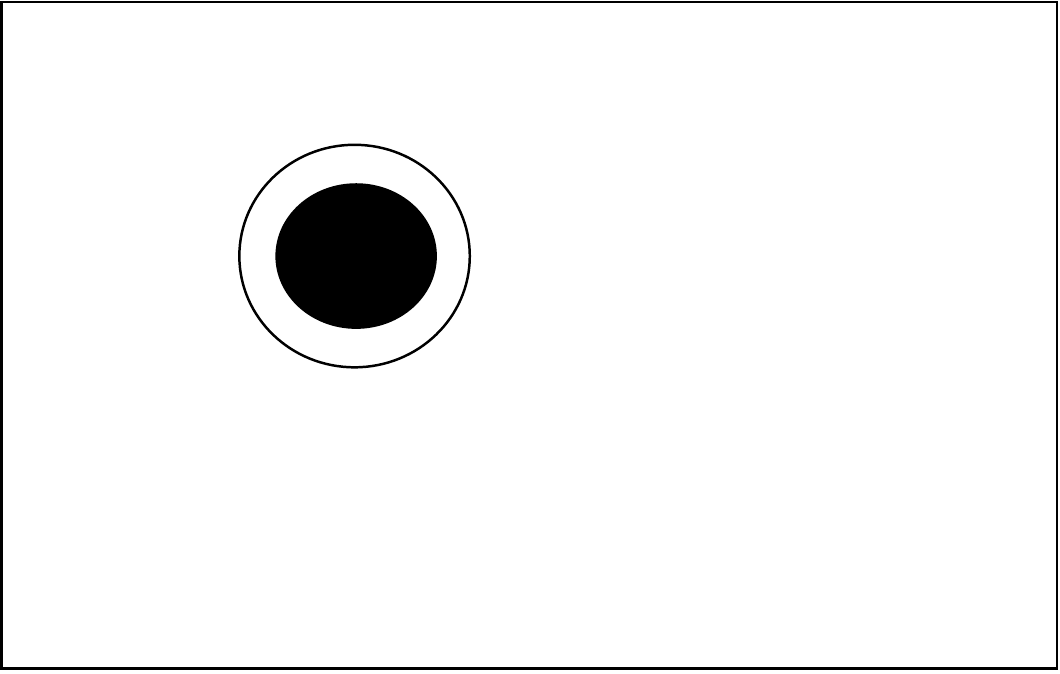}
		\centering {\small (b) ``bound''}
	\end{minipage}\hspace{20pt}
				\centering
	\begin{minipage}[t]{0.2\textwidth}
		\includegraphics[width=\textwidth]{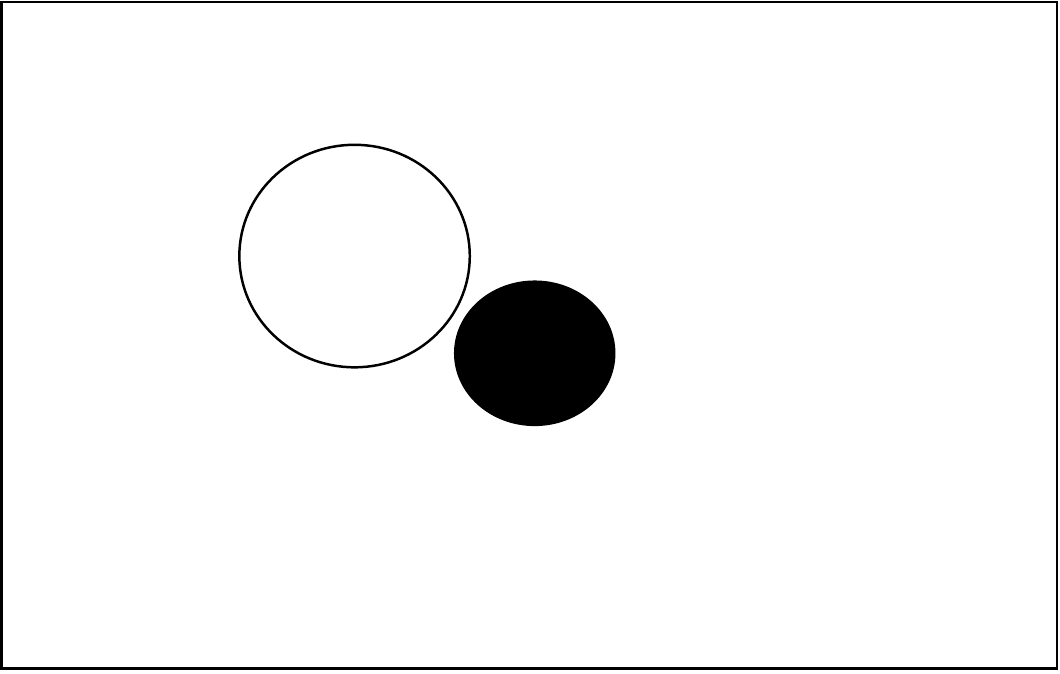}
		{\small (c) Spatial constellation after dissociation.}
	\end{minipage}
	\begin{minipage}[t]{0.2\textwidth}
		\includegraphics[width=\textwidth]{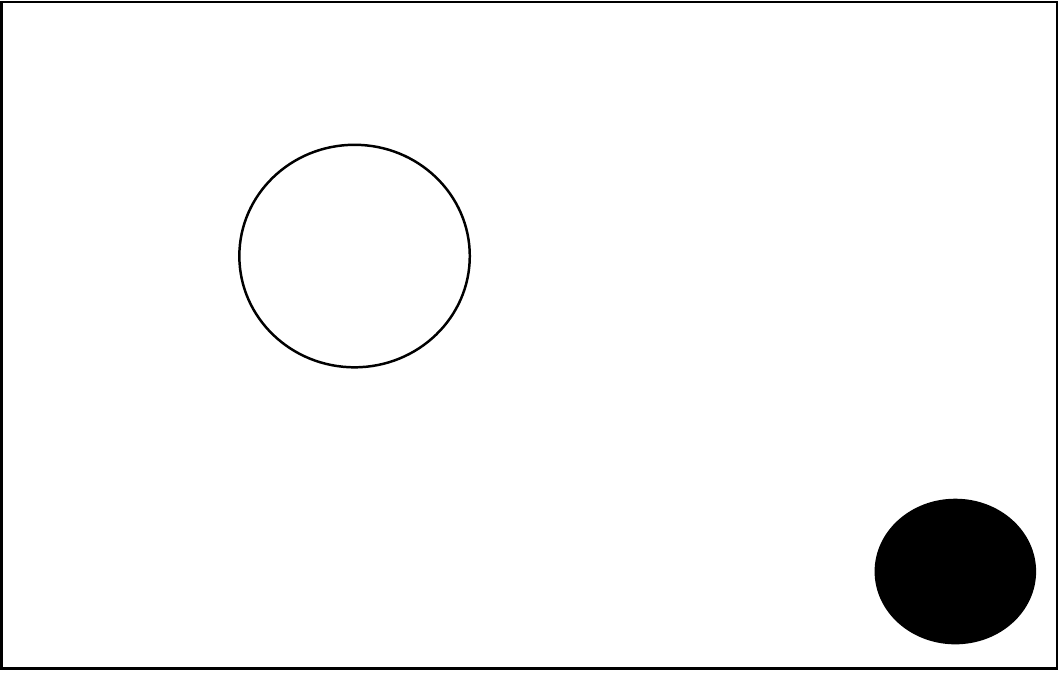}
		{\small (d) Spatial arrangement at arbitrary time.}
	\end{minipage}
	\captionof{figure}{(a) \& (b): Two possible macro states of a ligand-binding system. (c)\& (d): Rebinding effect: these two configurations represent the same macro state
		(``unbound'') and  are not distinguishable in model \eqref{eq:model_monovalent}, even though different binding probabilities are expected by the receptor-ligand-distance on the microscopic scale.}
	\label{fig:rebinding}
\end{minipage} 

The two possible macro states for a ligand-binding-system consisting of one receptor and one ligand are depicted in figure \ref{fig:rebinding} (a) \& (b).
We notice that the spatial arrangement of the receptor and the ligand in the unbound state is not included in the above model.
Therefore, we cannot distinguish if, at a given time, the receptor and the ligand are close to each other or not.
\\

By switching from the macroscopic to the microscopic point of view, we find out that the stochastic process modelled by \eqref{eq:model_monovalent} is actually \textbf{not} memoryless.
That is due to the spatial arrangement of the system after a receptor-ligand-complex dissociated.
Shortly after such a dissociation, it is more likely that the corresponding receptor and ligand will bind again, since they are still close to each other.
Such a binding shortly after a dissociation is called a \textbf{rebinding event}. The memory effect which thereby occurs is called \textbf{rebinding effect}.
On large timescales, this effect diminishes since the favorable spatial situation is not given anymore and the system is more likely to be rather mixed again.
Thus, Markovianity can be spoiled by the rebinding effect, as depicted in figure \ref{fig:rebinding} (c) \& (d). \\


In order to measure the magnitude of the rebinding effect, we \textbf{interpret} model \eqref{eq:model_monovalent} as the projection of a larger system. 
A crisp clustering does not yield a correct model and should be replaced by the fuzzy approach.
Accordingly, we consider the macro states ``unbound'' and ``bound'' as overlapping states.
This allows a micro state to be in the ``unbound'' macro state with a high degree of membership to the ``bound'' state, for instance shortly after a dissociation, which could be interpreted as an ``almost bound'' state.
Thus, if these states are strongly overlapping, then a high rebinding effect can be expected.
In the next sections, we quantify the rebinding effect by its relation to the magnitude of overlap of the conformations.
\\

The rebinding effect and its occurence in natural science has been described and analyzed by several authors\cite{goldstein1995approximating, vauquelin2010}.
In chemistry, it has been discussed in the context of clustered receptors and clustered ligands, e.g. multivalent systems\cite{care2011impact,fasting2012multivalency,von2016allosteric}.
A mathematical investigation of the rebinding effect has been realized by Weber et al\cite{weber2012,weber2014}.

\subsection{Measuring the Rebinding Effect}

We analyze the matrix representation $P_c = S^{-1} T$ of the Markov State Model.
The stochastic matrix $T$ represents the dynamical behaviour ot the process, though the Markov State Model differs from $T$ by
\begin{equation*}
S P_c(\tau) = T.
\end{equation*}
This ``deviation'' of the Markov State Model $P_c(\tau)$ from the coupling matrix $T$ is caused by the overlap of the membership functions, included in the matrix $S$.
If $S$ is equal to the identity matrix, then the Markov State Model is solely determined by $T$.
If $S$ is close to the identity matrix, then $P_c(\tau)$ is close to $T$ and not strongly influenced by $S$.
The more the overlap matrix $S$ differs from the identity matrix, the more the Markov State Model $P_c(\tau)$ differs from the transition matrix $T$.
This is due to the rebinding events.
The larger this deviation, the larger the occurring memory effects.
Thus, the rebinding effect, a memory effect provoked by a projection, can be measured by the matrix $S$.
The more the membership functions are overlapping, the more the matrix $S$ deviates from the identity matrix and thereby includes stronger memory effects.

Thus, the rebinding effect can be measured by the trace of the matrix $S$, being the sum of its diagonal elements. It can lie between $0$, implying very much rebinding, and $n$, implying no rebinding.
This approach to measure the rebinding effect has been introduced by Weber and Fackeldey\cite{weber2014} and will be used in the next section to detect a minimal bound for the rebinding effect included in a projected system.

\section{Optimization Problem}
\label{sec:optimization}

Commonly, we are mainly concerned to compute the projection of a large process and, of particular interest, to analyze how such a projection introduces memory effects in the clustered process.
In most of the cases though, we do not know the continuous transfer operator or infinitesimal generator describing a system. Instead, we are often given a finite matrix, for instance stemming from experimental data or as the solution of a differential equation.
In either case, such a finite matrix can be \textbf{interpreted} as a projection, since it is basically a model for an originally continuous process, describing the movement of molecules in $\mathbb{R}^3$.
\\

Assume we are in the situation that we only know the projected process $Q_c$. Nevertheless, we would like to know how much rebinding is included in that system, originating from the unknown projection.
Since we don't know on which membership functions the projection is based on, we can only compute an estimation for that. Considering all possible membership functions, how much rebinding is included \textbf{at least} in the system? In other words, how strongly overlapping are the membership functions at least?
\\

In \cite{weber2014} it is shown that the overlap matrix $S$ from \eqref{eq:matrix_representation} provides a measure for the quantity of the rebinding effect. In particular, being close to the identity matrix implies a low rebinding, while high outer diagonal elements of $S$ result in a high rebinding effect.
In order to reveal the actual impact of the rebinding effect, we set it in relation to the stability of the clustered system $Q_c$. 
Afterwards, we formulate an optimization problem in order to deduce a lower bound for the rebinding effect included in a given system.
For the sake of simplicity, we assume in the further course that the transition rates can be measured experimentally.
Accordingly, we examine the given transition rate matrix $Q_c$ of a process.

\subsection{Relevance of the Rebinding Effect}

If the eigenvalues $\xi_i$ of $Q_c$ are close to $0$, then the macro states are very stable in the sense that the probability to stay inside of such a state is close to $1$.
The trace of $Q_c$ corresponds to the sum of the dominant eigenvalues of $Q$.
Thus, we can measure the \textit{stability} of the molecular system by the quantity $F := - \mathrm{trace}(Q_c)$.
If $F$ is close to $0$, then the system is very stable, while it is less stable for a high value of $F$.
We want to set the stability $F$ in relation to the measure of the rebinding effect, the overlap matrix $S$.

Let $Q_c$ be the projected infinitesimal generator of a process and $P_c(\tau)$ the corresponding projected transfer operator with the matrix representation $P_c(\tau) = S^{-1}T$, then according to \cite{weber2014}, the quantity $F := - \mathrm{trace}(Q_c)$ can be measured by
\begin{equation}
\label{eq:stability}
F = \tau^{-1}(\mathrm{log}(\mathrm{det}(S)) - \mathrm{log}(\mathrm{det}(T))).
\end{equation}
The coupling matrix $T$ describes the stochastic movement of the process and in particular, encodes the metastable behaviour between the conformations. Large diagonal elements result in a strong metastability and a slow process, while higher outer diagonal elements lead to faster transitions between the metastable sets.
On the other hand, the overlap matrix $S$ merely includes informations about the crispness of the membership functions, implying the magnitude of the rebinding effect.

Equation \eqref{eq:stability} shows that \textbf{both} determinants of $S$ and $T$  influence the stability of the system, though in opposite directions.
If $\mathrm{det}(T)$ is close to $1$, then $F$ is low and consequently the process is rather stable.
If $\mathrm{det}(T)$ is small, then the process is rather unstable, since $F$ is high.
A high determinant of $T$ leads to a high metastability of the system and thus describes a slower process, while a low determinant implies higher outer diagonal elements of $T$ and thus, makes the process faster. 

In contrast, if $\mathrm{det}(S)$ is close to $1$, then the first term in \eqref{eq:stability} vanishes and hence, $S$ barely contributes to the stability, which is instead mainly determined by $T$.
On the other hand, if $\mathrm{det}(S)$ is close to $0$, the system becomes more stable.
This means that a higher overlap of the membership functions, and thus a \textbf{strong rebinding effect}, leads to a more stable process.
\\

At first sight, it sounds plausible to equalize the stability of a system to its slowness.
A slow system has rare transitions and thereby implies a stable system.
However, a stable system does \textbf{not} necessarily imply a slow system.
Instead, a rather fast system can gain a certain stability by the rebinding effect.
The ``fast'' system has frequent transitions between its metastable sets. However, in case of a strong rebinding, the quitting of a metastable set can with high probability be followed by an immediate return to the previous state.
Thus, the rapidness of the process can to a certain extent be compensated by the rebinding effect.
Concluding, we can differentiate between two factors leading to a high stability:
\begin{itemize}
	\item $\mathrm{det}(T)$ high:
	The conformations have a high metastability and are well-separated. Therefore, transitions between the metastable sets are rare and the process is slow.
	\item $\mathrm{det}(S)$ low: A high rebinding effect makes the process more stable, since transitions out of a metastable set can be compensated by a fast transition back. In particular, a rapidly mixing process, $\mathrm{det}(T) \ll 1$, can be stabilized by the rebinding effect.
\end{itemize}
A stable system is naturally reached by a strongly metastable matrix $T$, though can likewise be obtained for a weaker metastable matrix $T$, if much rebinding is included.

\subsection{Lower Bound for the Rebinding Effect}

In order to determine the stability of a system, it is of interest to know how much rebinding is included.
We compute a lower bound to find out how much rebinding we are guaranteed \textbf{at least}.
In order to derive an optimization problem, let us first remember how $S$ is determined.
The transition matrix $P \in \mathbb{R}^{m \times m}$ is projected onto a finite-dimensional state space via membership functions $\chi$ as a linear combination of the dominant Schur vectors with a regular matrix $A$. The choice of the matrix $A$ determines $S$ and in particular the magnitude of rebinding.
In order to estimate the rebinding effect included in a system, we take into consideration all \textit{feasible} transformation matrices $A$, see\cite{weber2006meshless}.

Similar to\cite{weber2014}, we formulate an optimization problem to reveal which choice of $A$ results in the \textbf{lowest} rebinding effect, measured by an \textit{optimal matrix} $\Sopt$.
This problem is equivalent to finding the largest possible determinant of $S$.
\\

We are interested in the rebinding effect included in the clustered system $Q_c$.
If we know the employed membership functions $\chi$ or the transformation matrix $A$, then we can easily compute the \textbf{real} rebinding effect which is encoded in the overlap matrix $S = D^{-1} A^T A$.
\\
Given a finite matrix $Q_c$ and a Schur decomposition $\Xi$ with the corresponding Schur vectors as columns of the matrix $X$, the starting point to construct the optimization problem is given by
\begin{equation}
\label{eq:eigenvalue_problem}
Q_c X = X \Xi.
\end{equation}
Then $\Xi$ is of block-triagonal shape, according to section \ref{sec:projection}.
Since we assume that the dominant eigenvalue $\lambda_1 = 1$ is unique, 
the first column of $X$ corresponds to the first Schur vector $X_1 := (1,\dots, 1)^T$.
By \eqref{eq:matrix_representation4}, we see that $A^{-1}$ is a matrix of Schur vectors for the Schur decomposition $\Xi$ as well.
\\

Assuming a reversible process, then the Schur decomposition is equal to the spectral decomposition and results in a \textbf{diagonal} eigenvalue matrix $\Xi$.
Therefore the columns of $A^{-1}$ consist of multiples of the eigenvectors $X_j$, yielding
\begin{equation}
\label{eq:Ainv_eigen}
A^{-1} =
\begin{pmatrix}
1 	  & & & \\
\vdots & \alpha_2 X_2 & \cdots & \alpha_n X_n \\
1	  & & &
\end{pmatrix}
\end{equation}
with $\alpha_1 = 1$ and $\alpha_2,\dots,\alpha_n \in \mathbb{R}$.
However, the diagonal shape of $\Xi$ is \textbf{not} guaranteed  for a non-reversible process. Instead, it may contain $2 \times 2$-blocks, which have to be considered. For the case of such a $2 \times 2$-block, the two associated Schur vectors are not linear independent and thus cannot be simply built as a multiple of the corresponding Schur vectors from $X$.
\\
We therefore take a different path, by employing an optimization procedure.
To keep things simple, we consider a Schur decomposition with Schur Matrix $\Xi \in \mathbb{R}^{3\times 3}$
\begin{equation}
\label{eq:schur_structure}
\Xi = 
\begin{tikzpicture}[baseline=-0.5ex]
\matrix [matrix of math nodes,left delimiter=(,right delimiter=)] (m)
{
	0    & 0  & 0     \\
	0 & *   &  \spadesuit     \\
	0   & \clubsuit    & *       \\
} ;

\draw (m-2-2.north west) rectangle (m-3-3.south east);
\end{tikzpicture},
\end{equation}
where the $*$ denote non-zero entries (the real parts of the eigenvalues of $Q_c$) and $\clubsuit$ is non-zero only in the case of complex eigenvalues. $\clubsuit$ and $\spadesuit$ equal zero in the reversible case. Note, that the first row of $\Xi$ is always zero in Markov processes.
Given the matrix $X$ of Schur vectors associated to $\Xi$, we aim to reveal the necessary structure of $A^{-1}$ such that 
\[
A Q_c A^{-1} = \Xi.
\]
Because the membership vectors $\chi$ should sum up to one, the first column of $A^{-1}$ consists only of ones, i.e. 
\[
A^{-1}=\begin{tikzpicture}[baseline=-0.5ex]
\matrix [matrix of math nodes,left delimiter=(,right delimiter=)] (m)
{
	1    & *  & *     \\
	1 & *   & *      \\
	1   & *    & *       \\
} ;

\end{tikzpicture}.
\]
The fact, that $S$ is a stochastic matrix is
leads to a further side constraint
\begin{equation}
\label{eq:opti_schur3}
\mbox{
	\boxed{ S_{ij} \geq 0 \ \ \forall i,j}.
}
\end{equation}
This means, that $(A^TA)_{i,j} \geq 0\, \forall i,j$.

Let $e_1=(1\, 0\, 0)^T$ and $e=(1\, 1\, 1)^T$ we can summarize these constraints in the following set
\[
{\cal C}=\{ A \in \mathbb{R}^{3\times 3} | A Q_c A^{-1} = \Xi \text{ and } A^{-1}e_1=e \text{ and } A^TA \geq 0 \}.
\]

Based on these relations, we can formulate an optimization problem. We know that a determinant of $S$ close to $1$ results in a low rebinding effect. Thus, in order to find a lower bound, we try to maximize $\mathrm{det}(S)$, or equivalently minimize $|\mathrm{det}(S)-1|$, since $S$ is a stochastic matrix having $1$ as largest possible determinant.
Then the \textit{objective function} of the optimization problem is given by
\begin{equation}
\label{eq:optimization}
\mbox{
	\boxed{ \mathrm{min}_{A \in {\cal C}} |\mathrm{det}(S) -1|},
}
\end{equation}
where several \textit{side constraints} have to be included, which lead to a stochastic matrix $S$.

A \textit{feasible solution} of this optimization problem is a matrix $S$ fullfilling all side contraints, but not necessarily being an optimum.
Any feasible solution of optimization problem \eqref{eq:optimization} will be called a \textit{real overlap matrix} $\Sreal$, while an actual optimum will be called an \textit{optimal overlap matrix} $\Sopt$. Clearly, we get $\mathrm{det}(\Sreal) \leq \mathrm{det}(\Sopt) \leq 1$.

\subsection{Interpretation}

The real rebinding effect is high if the determinant of $\Sreal$ is low. Thus, a small determinant of $\Sopt$ implies a high rebinding effect, while a large determinant of $\Sopt$ gives us only few information about the actual quantity of the rebinding effect, it could be either large or small.
Unfortunately, a reversible process $Q_c$ yields a trivial solution of optimization problem \eqref{eq:optimization} and therefore, provides us with no information, as it has been shown in{\cite{weber2014}}.
That means that for every such process, it is possible to find a transformation matrix $A$ which causes no rebinding.
Consequently, a nontrivial estimation for the rebinding effect can be obtained only for a nonreversible system $Q_c$.
In particular, only systems with at least three states are of interest to examine, since $Q_c$ is reversible for $n=2$.
For instance, the  example from section \ref{sec:rebinding} describing a receptor-ligand system on two macro states ``bound'' and ``unbound'' yields the trivial solution.

Optimization problem \eqref{eq:optimization} is a generalized version of the minimization problem for reversible processes from\cite{weber2014}. Due to the block-triagonal shape of $\Xi$, it requires a case distinction of the different Schur blocks. However, it includes reversible as well as non-reversible processes. For a reversible system, where the Schur decomposition consists of $1 \times 1$-blocks corresponding to the dominant eigenvalues, they coincide.
Thus, with this generalization, we can compute the minimal rebinding effect for \textbf{any} system, independent of the reversibility or non-reversibility of the original process.
The quality of this estimation  will be evaluated in the next chapter by means of an exemplary reversible process, which will be slightly perturbed to non-reversibility by introducing such a $2 \times 2$-block.
The solution of optimization problem \eqref{eq:optimization} will be computed in the following for some illustrative examples.

\section{Numerical Examples}
\label{sec:numerical}

The results from section \ref{sec:optimization} will be verified analyzing two illustrative examples: an artificial process indicating the role of the non-reversibility towards the minimal rebinding effect and a `real-world' process describing a chemical reaction.

\subsection{Artificial Example}\label{sec:art}

We consider a system given by the Schur decomposition
\begin{equation}
\label{eq:schur_decomp}
\Lambda =
\begin{pmatrix}
1 & 0 & 0 & 0 & 0 \\
0 & 0.99 & \epsilon & 0 & 0 \\
0 & -\gamma & 0.98+\delta & 0 & 0 \\
0 & 0 & 0 & 0.005 & 0 \\
0 & 0 & 0 & 0 & 0.001
\end{pmatrix},
\end{equation}
with $\epsilon, \gamma, \delta > 0$. The corresponding transition matrix is computed by $P = X \Lambda X^{-1}$, with a set of Schur vectors $X$. If $\Lambda$ is a diagonal matrix, then this equation represents the eigenvalue problem of a reversible process $P$. By introducing non-zero values for $\epsilon, \gamma$ and $\delta$, the system gets non-reversible.
This example is of particular interest, since we examine different systems, yet having the same Schur vectors and very similar Schur decompositions. However, these small changes in the Schur decomposition lead to different results when it comes to computing the minimal rebinding effect.

Having three dominant eigenvalues, the matrix \eqref{eq:schur_decomp} describes a system on three metastable sets. Accordingly, we examine different clustering on a  three-dimensional state space. For that aim, we employ several transformation matrices $A \in \mathbb{R}^{3 \times 3}$, turning the dominant Schur vectors $X \in \mathbb{R}^{5 \times 3}$ into membership functions $\chi \in \mathbb{R}^{5 \times 3}$.
We generate $200$ random feasible transformation matrices $A$ and examine the rebinding effect caused by this projection. Finally, we compare this real rebinding effect with the minimal rebinding effect included in the clustered system, as the solution of optimization problem \eqref{eq:optimization}. \\

\begin{minipage}{\linewidth}
	\centering
	\begin{minipage}[t]{0.3\textwidth}
		\includegraphics[width=\textwidth]{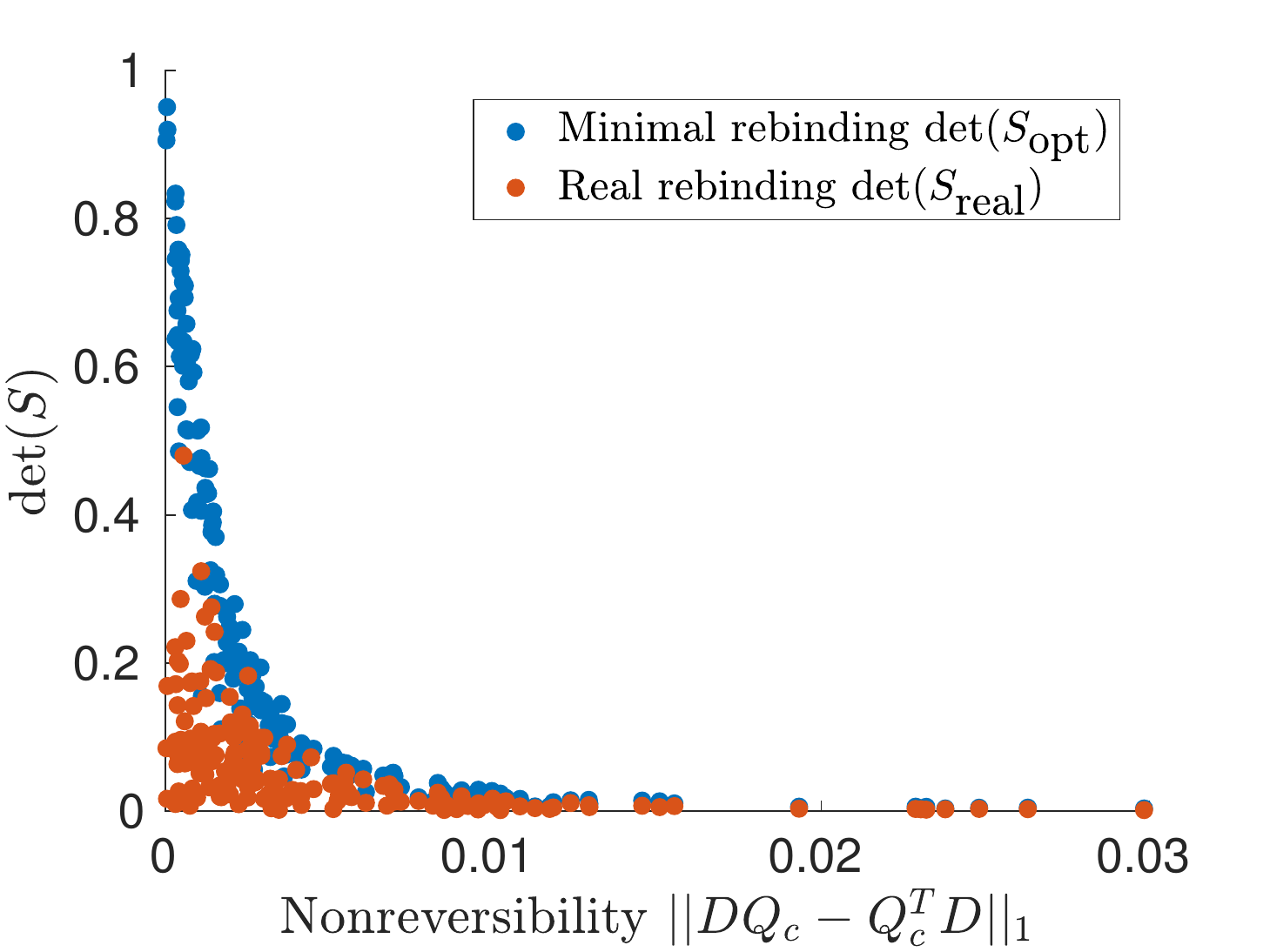}
		(a) Minimal and real rebinding effect compared to the degree of non-reversibility of $Q_c$.
	\end{minipage}
	\hspace{20pt}
	\begin{minipage}[t]{0.3\textwidth}
		\includegraphics[width=\textwidth]{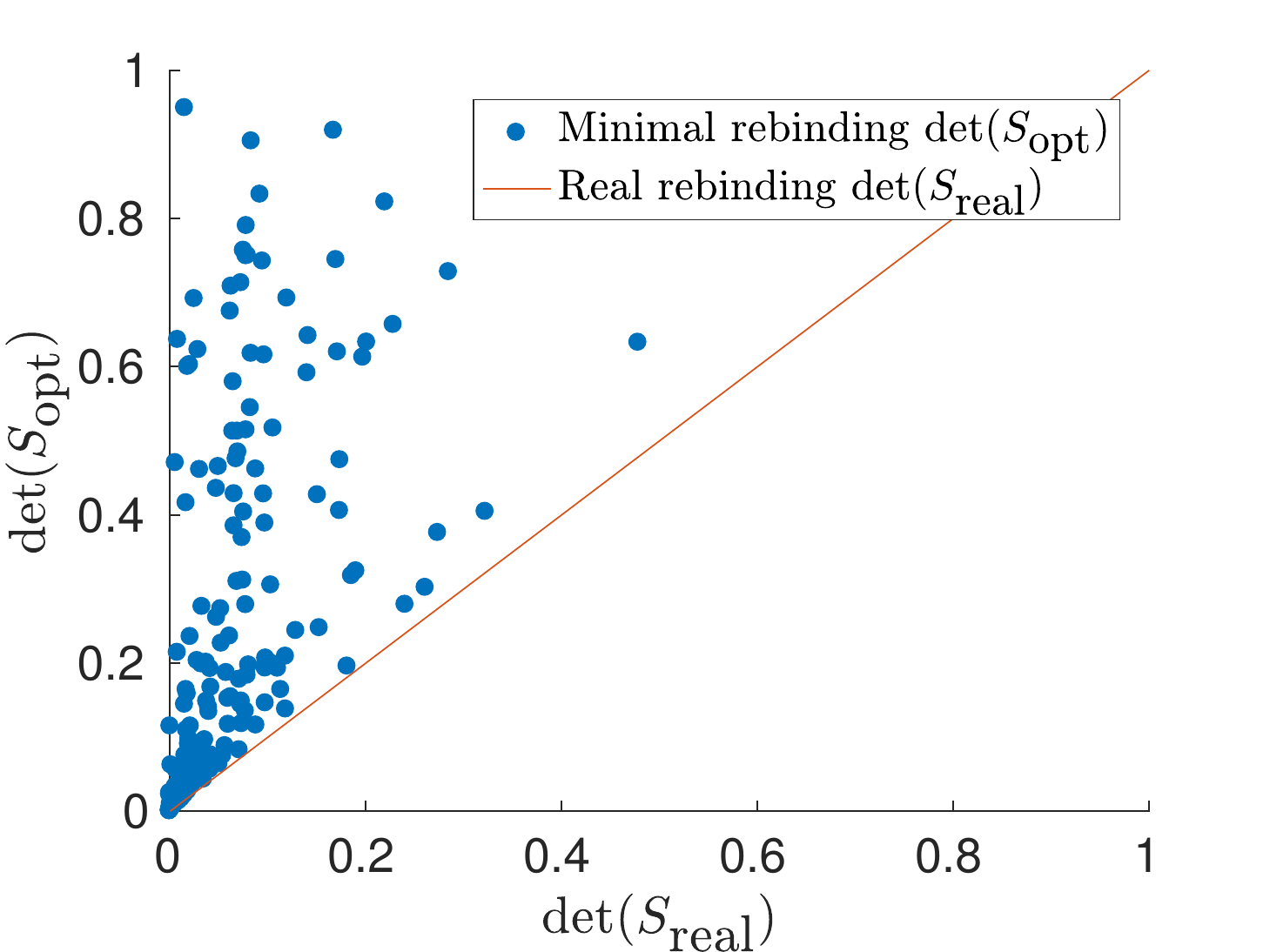}
		(b) The minimal rebinding effect compared to the real rebinding effect included in $Q_c$.
	\end{minipage}
	\captionof{figure}{The system $P$ is clustered with $200$ randomly generated transformation matrices $A$ for the parameters $\epsilon = \delta = \gamma = 0$.}
	\label{fig:reb_example}
\end{minipage} \\

The degree of non-reversibility can be measured by $\Vert D Q_c - Q_c^T D \Vert_1$.
We investigate the minimal rebinding effect $\mathrm{det}(\Sopt)$ depending on the degree of non-reversibility of the clustered system $Q_c$.
For all examined systems, we observe a considerable correlation between the lower bound of the rebinding effect and the non-reversibility of the system: the more non-reversible the system, the larger the minimal rebinding effect.
For the reversible case, where all outer-diagonal elements in \eqref{eq:schur_decomp} are $0$, this correlation is very strong, see figure \ref{fig:reb_example} (a), though it can be a rather good or a rather bad estimation, see figure \ref{fig:reb_example} (b).

Inserting a small outer-diagonal perturbation $\epsilon = 0.004$ leads to a non-reversible process. For different clusterings, the minimal rebinding effect behaves similar to the reversible case, yet being slightly more uneven, see figure \ref{fig:reb_example_eps}. \\

\begin{minipage}{\linewidth}
	\centering
	\begin{minipage}[t]{0.3\textwidth}
		\includegraphics[width=\textwidth]{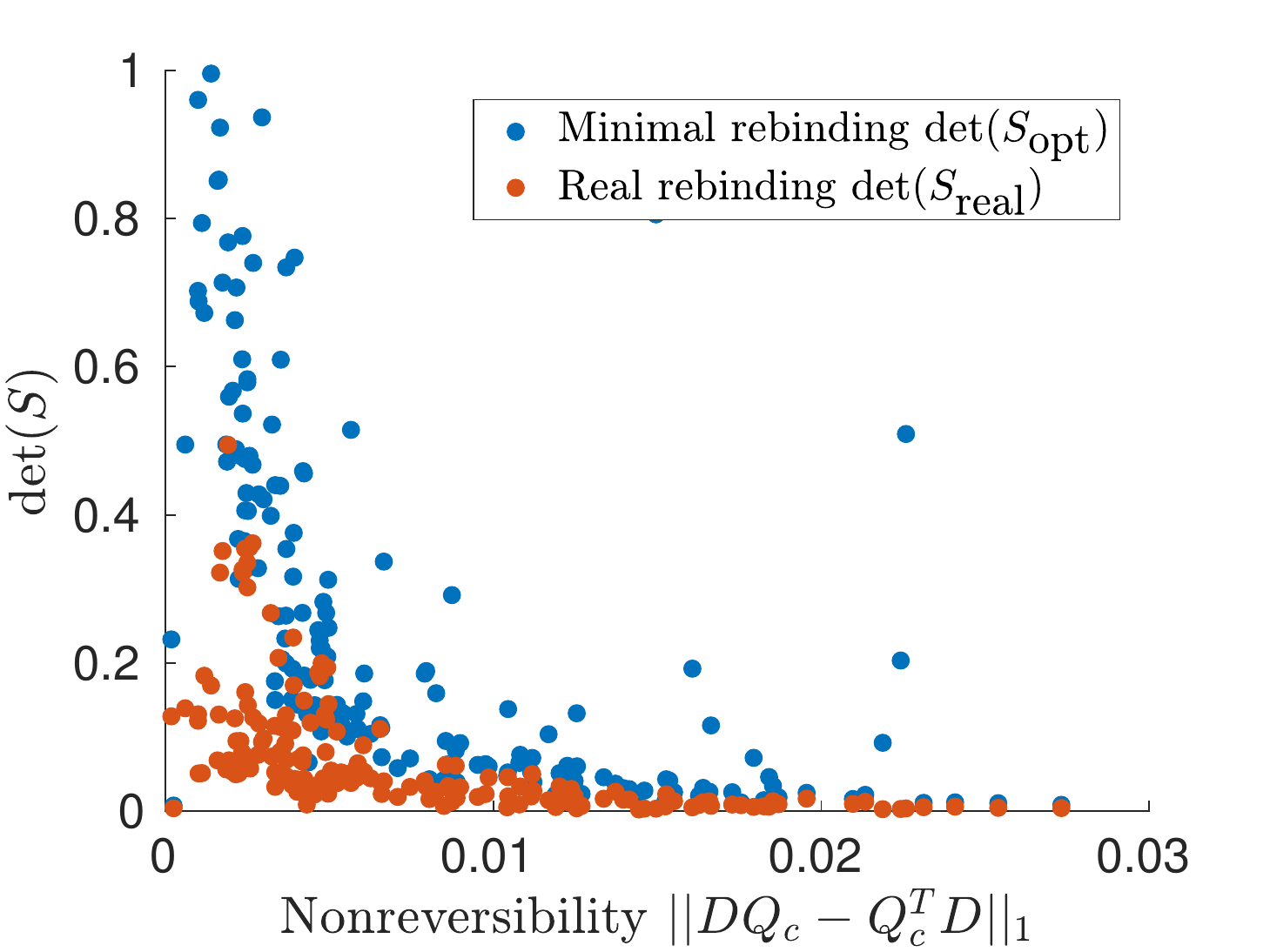}
		(a) Minimal and real rebinding effect compared to the degree of non-reversibility of $Q_c$.
	\end{minipage}
	\hspace{20pt}
	\begin{minipage}[t]{0.3\textwidth}
		\includegraphics[width=\textwidth]{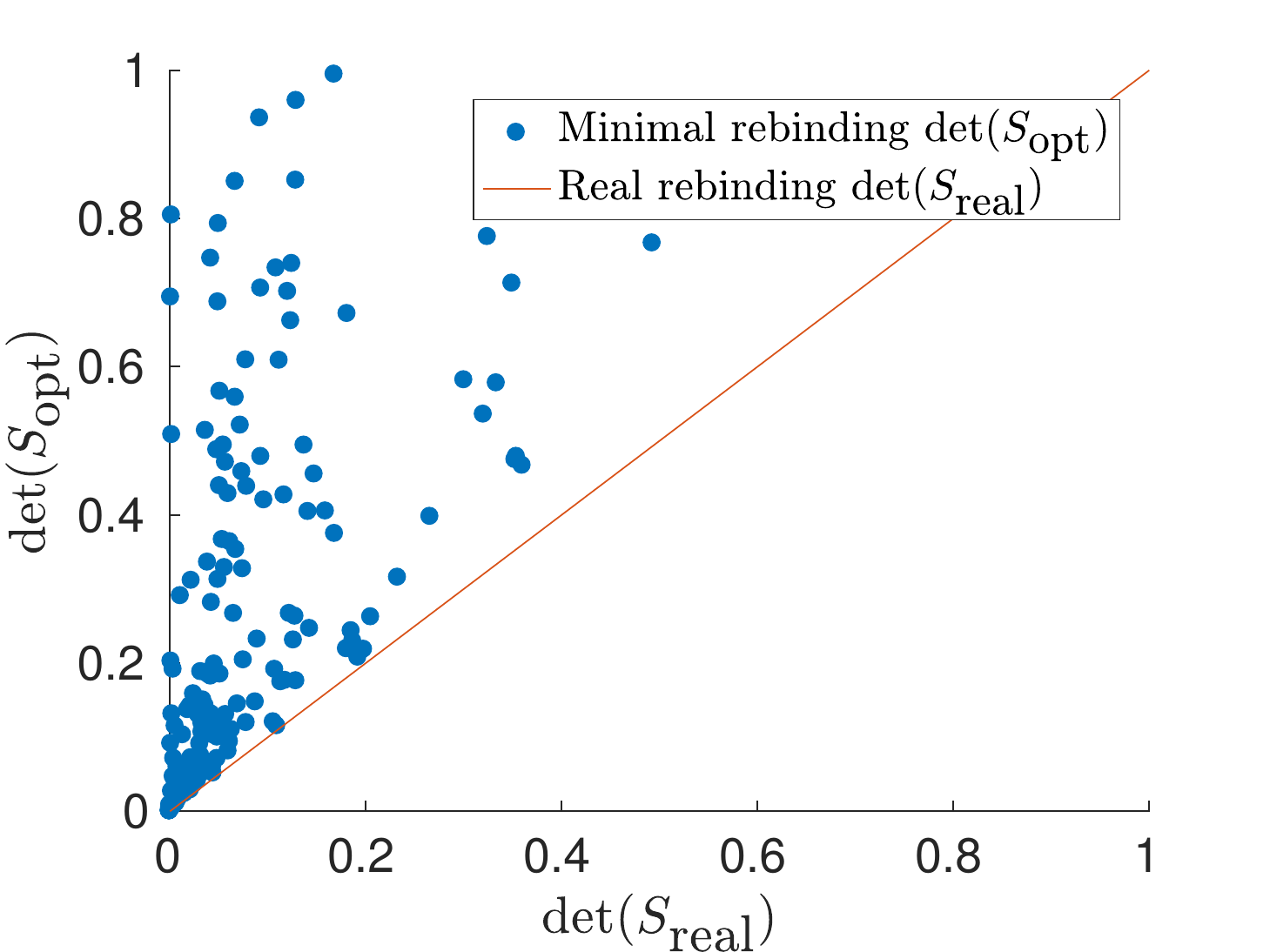}
		(b) The minimal rebinding effect compared to the real rebinding effect included in $Q_c$.
	\end{minipage}
	\captionof{figure}{The system $P$ is clustered with $200$ randomly generated transformation matrices $A$ for the parameters $\epsilon = 0.004$, $\delta = \gamma = 0$.}
	\label{fig:reb_example_eps}
\end{minipage} \\

A further perturbation $\delta = 0.01$ increases the non-reversibility of the system and leads to a non-diagonalizable matrix $P$. The minimal rebinding effect for the clusterings is presented in figure \ref{fig:reb_example_delta}. \\

\begin{minipage}{\linewidth}
	\centering
	\begin{minipage}[t]{0.3\textwidth}
		\includegraphics[width=\textwidth]{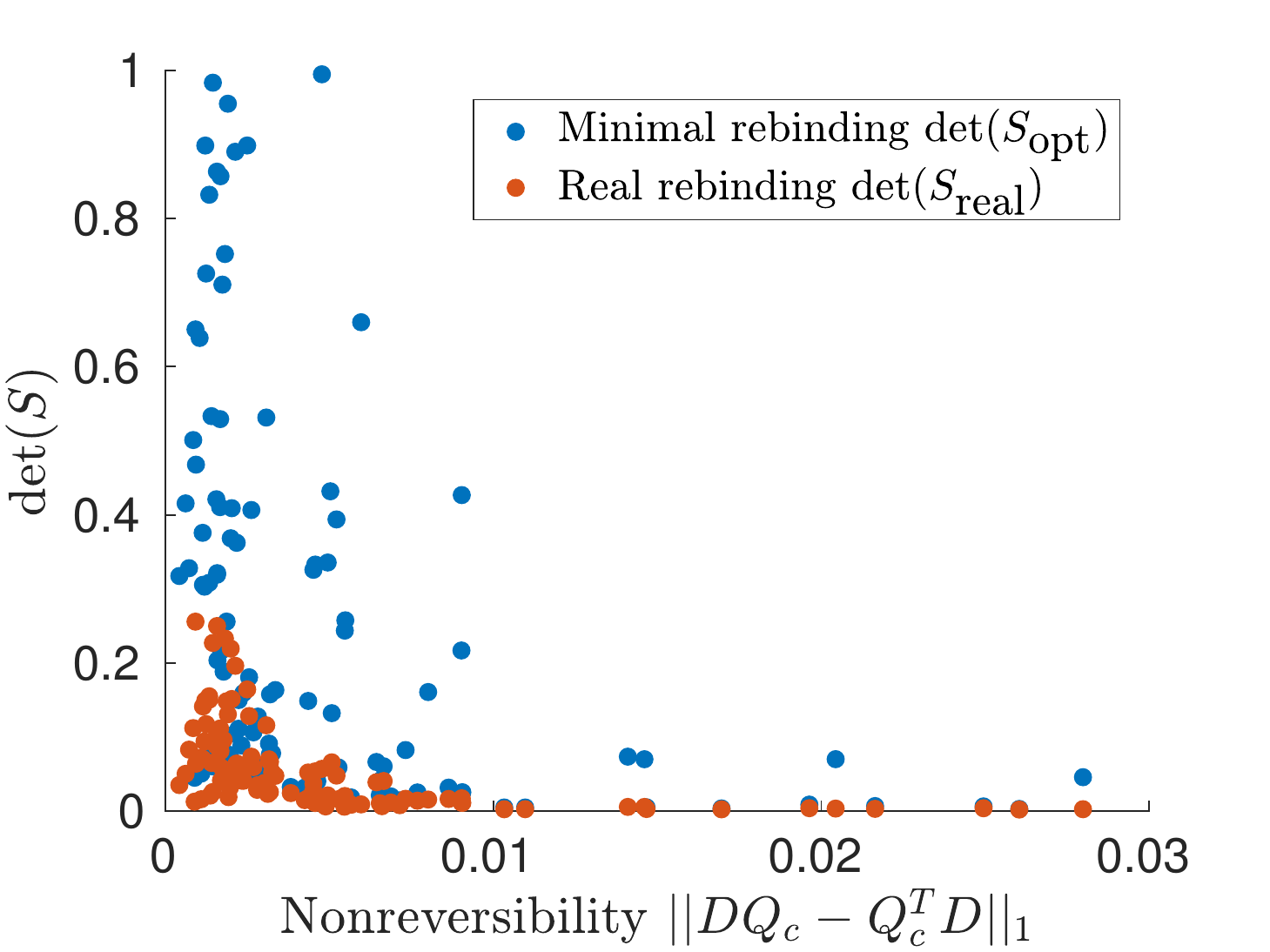}
		(a) Minimal and real rebinding effect compared to the degree of non-reversibility of $Q_c$.
	\end{minipage}
	\hspace{20pt}
	\begin{minipage}[t]{0.3\textwidth}
		\includegraphics[width=\textwidth]{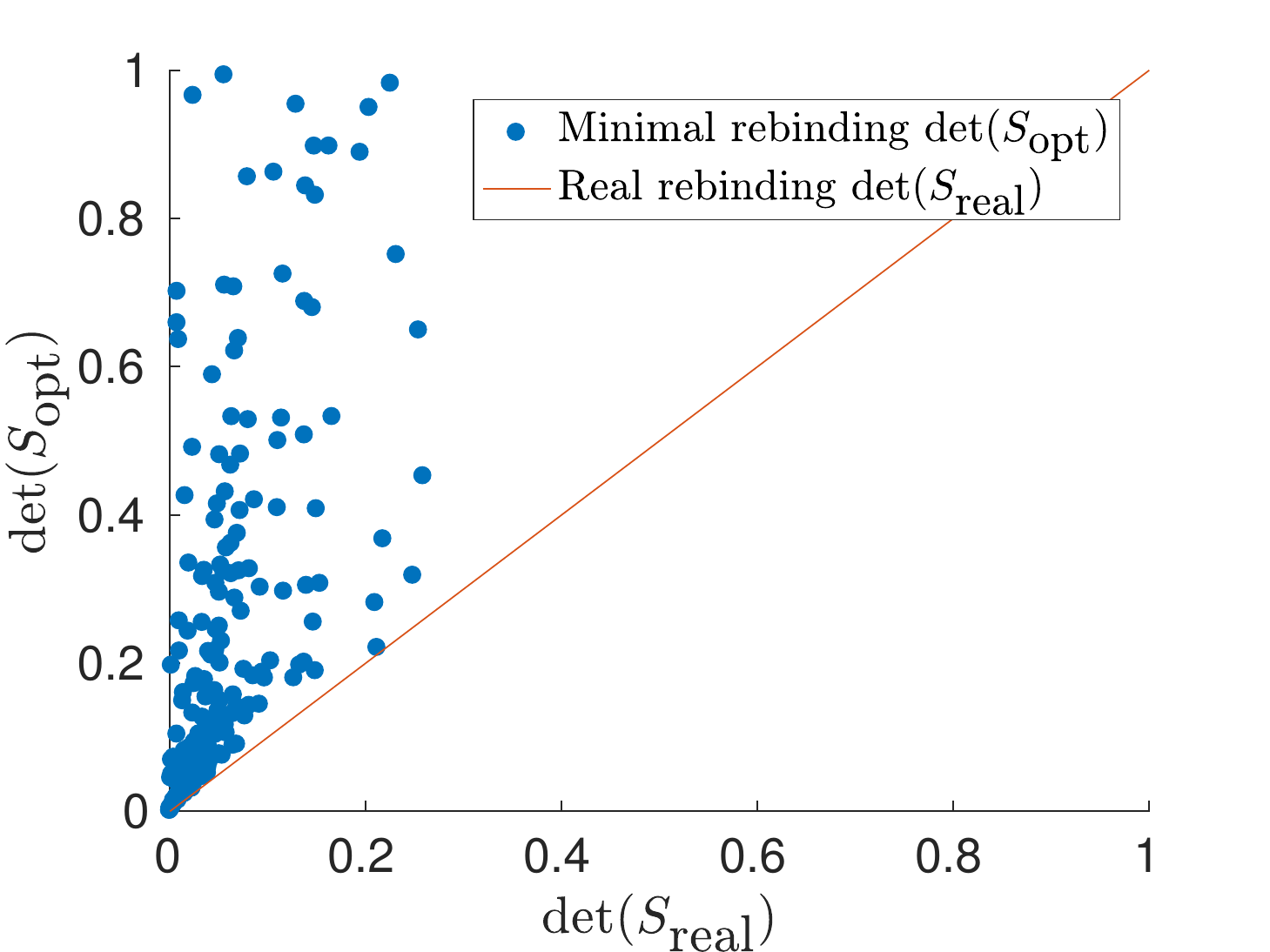}
		(b) Minimal rebinding effect compared to the real rebinding effect included in $Q_c$.
	\end{minipage}
	\captionof{figure}{The system $P$ is clustered with $200$ randomly generated transformation matrices $A$ for the parameters $\epsilon = 0.004$, $\delta = 0.01$, $\gamma = 0$.}
	\label{fig:reb_example_delta}
\end{minipage} \\

The general tendency of the results is similar for all tested parameters:
while the quality of the estimation can be either good or bad, there is a clearly visible correlation between the minimal rebinding effect $\mathrm{det}(\Sopt)$ and the non-reversibility of $Q_c$. However, this correlation seems to diminish the more we `perturb' the original process from reversibility. This weakened correlation implies that for originally non-reversible systems, the quality of the estimation is less predictable.

\subsection{Electron Densities}\label{sec:ED}

The occurrence of some kind of rebinding effect can be observed in all different types of processes when projecting them.
The actual meaning of this effect has to be interpreted for each system individually.
We present a process describing the change of electron densities during a pericyclic chemical reaction, examined in \cite{Schild2013,weber2017coherent}.

\begin{figure}[ht!]
	\centering
	\includegraphics[width=0.4\textwidth]{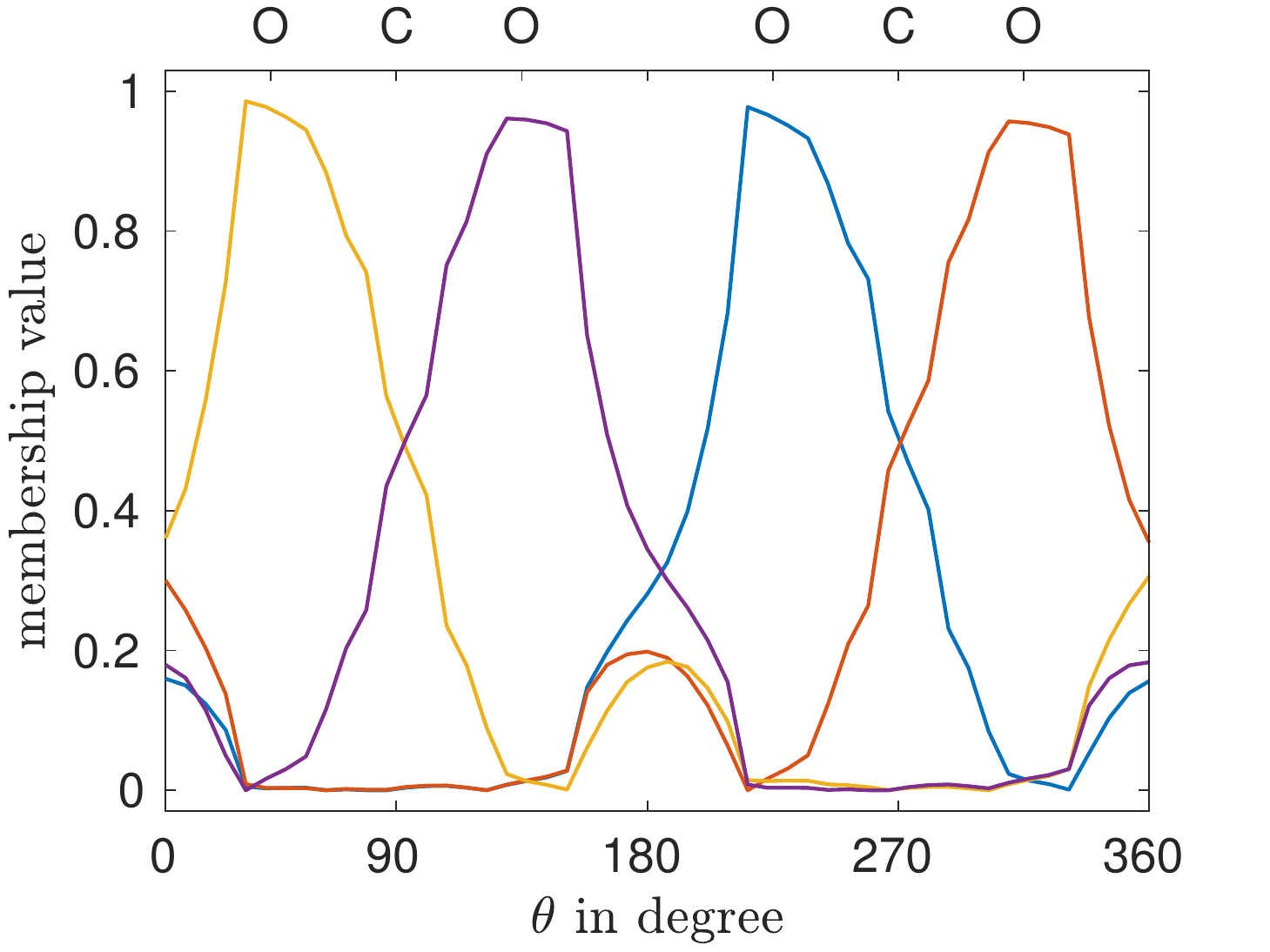}
	\caption{Membership functions obtained by GenPCCA.}
	\label{fig:electron_membership}
\end{figure}

Formic acid is a molecule consisting of one carbon atom C, two oxygen atoms O and two hydrogen atoms H.
\begin{figure}[!ht]
	\centering
	\includegraphics[width=0.29\textwidth]{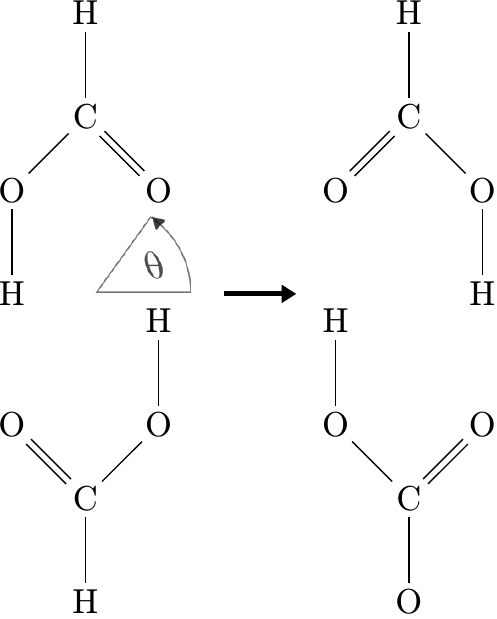}
	\caption{Chemical reaction in formic acid dimer.}
	\label{fig:acid_dimer}
\end{figure}
In such a system, reactions between the individual molecules take place, building hydrogen-bonded dimers, as depicted in figure \ref{fig:acid_dimer}. An H-atom which is attached to an O-atom moves to the O-atom of another molecule and vice versa.
These reactions are caused by double proton tunneling\cite{schild2013}.
During that process, the electron density changes accordingly.
The formic acid dimer cannot satisfactorily be described by one single 
Lewis formula. The two forms presented in figure ~\ref{fig:acid_dimer} are 
mesomeric formulas of this dimer. Thus, it is expected that the 
separation between these two types cannot be strict and the rebinding 
effect should be relevant.
This process can be represented by a reversible transition matrix $P$ consisting of the time-dependent electron densities $\pi(t)$, as described by\cite{weber2017coherent}.
Clustering it into four metastable sets using GenPCCA and transforming it into a transition rate matrix yields
\begin{equation*}
Q_c = 
\begin{pmatrix}
-2.0040  &  1.6859  &  0.1490  &  0.1690 \\
1.6192 &  -2.0010  &  0.1724  &  0.2095  \\
0.1451  &  0.1747 &  -1.9548  &  1.6350  \\
0.1632  &  0.2106  &  1.6217  & -1.9955
\end{pmatrix}.
\end{equation*}
The membership functions of this clustering are represented depending on the angle $\theta$ in figure \ref{fig:electron_membership}.
We notice that the four metastable conformations correspond to the angular regions of the O-atoms.
That means that high electron densities are detected around the O-atoms, which is plausible since the H-atoms tend to be attached to an O-atom.
Even though clustered with GenPCCA, having the objective of maximizing the crispness, we identify rather strongly overlapping membership functions in figure \ref{fig:electron_membership} and expect a high rebinding effect.
However, solving optimization problem $\eqref{eq:optimization}$ for $Q_c$ yields a lower bound
\begin{equation*}
\mathrm{det}(\Sopt) = 1,
\end{equation*}
providing us with no information, which can be explained by the reversibility of the clustered system, observed by $\Vert DQ_c - Q_c^T D \Vert_1 = 0$.
Knowing the membership functions $\chi$ and the stationary distribution $\pi$ of the original process, we can compute the \textbf{real} rebinding effect as
\begin{equation*}
\mathrm{det}(\Sreal) = \mathrm{det}(D^{-1} \langle \chi, \chi \rangle_\pi) = 0.2925,
\end{equation*}
corresponding to a strong overlap of the membership functions. 
Rebinding in this context can be interpreted similar to the rebinding in receptor-ligand-systems:
Shortly after a $H$-atom unbinds from an O-atom moving forward to the O-atom of a different molecule, it is still spatially close and attracted to its previous O-atom and therefore can rebind to it.
That is one factor contributing to the stability of the four conformations.
The quantitative influence of the rebinding effect on the stability of the clustered system is visualized in figure \ref{fig:electron_stability1} and \ref{fig:electron_stability2} for two different lag-times $\tau_1 = 0.2$ and $\tau_2 = 0.001$.
The metastability of the coupling matrix $T$ is enhanced by the significant overlap of the membership functions, yielding a strongly metastable transition matrix $P_c = S^{-1} T$. 
This confirms the result from section \ref{sec:rebinding}: the rebinding effect stabilizes a system by ``compensating'' a rather weak metastability of the conformations. \\

\begin{minipage}{\linewidth}
	\centering
	\begin{minipage}[t]{0.29\textwidth}
		\includegraphics[width=\textwidth]{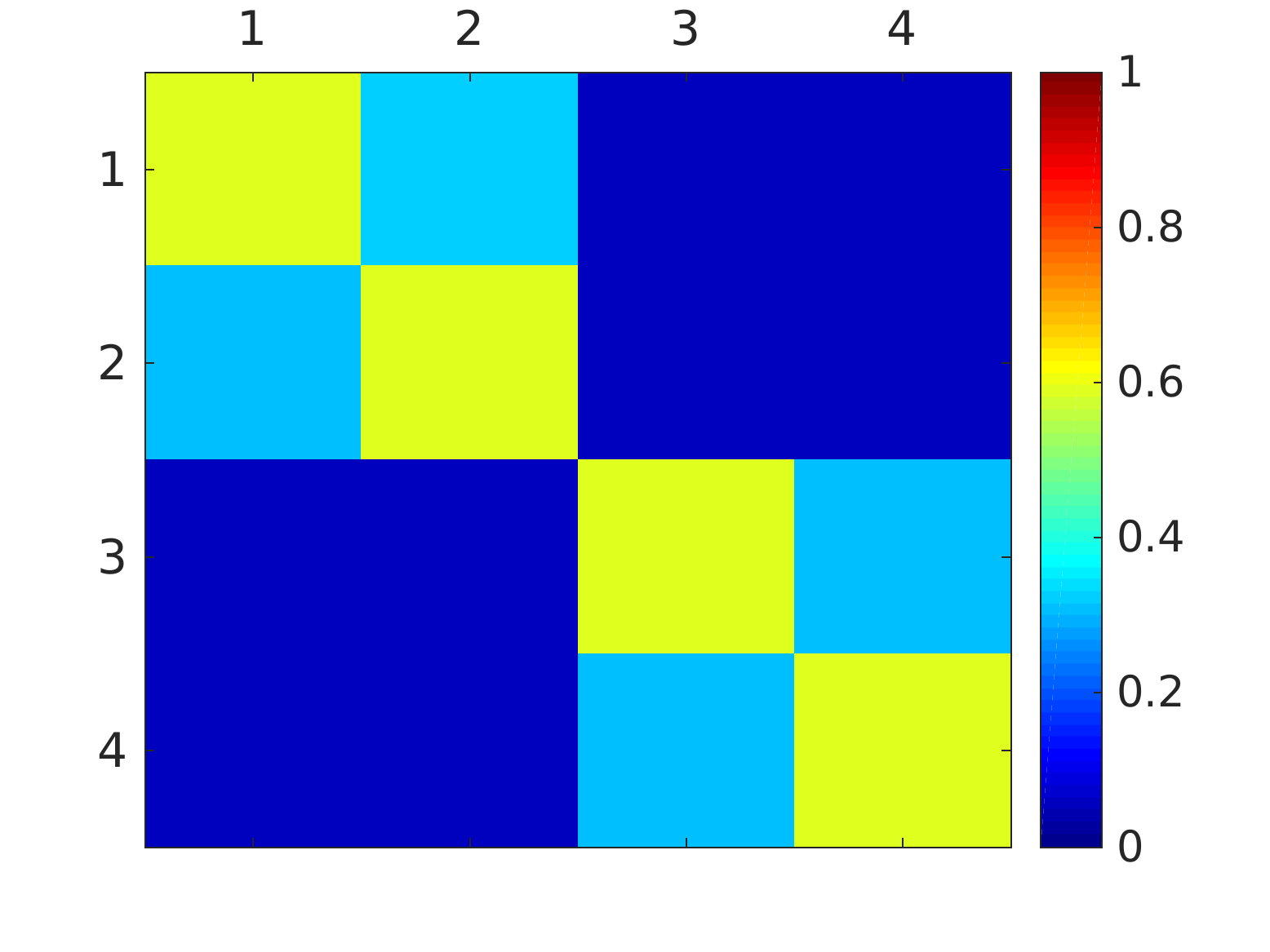}
		(a) Weakly metastable coupling matrix $T$.
	\end{minipage}
	\hspace{20pt}
	\begin{minipage}[t]{0.29\textwidth}
		\includegraphics[width=\textwidth]{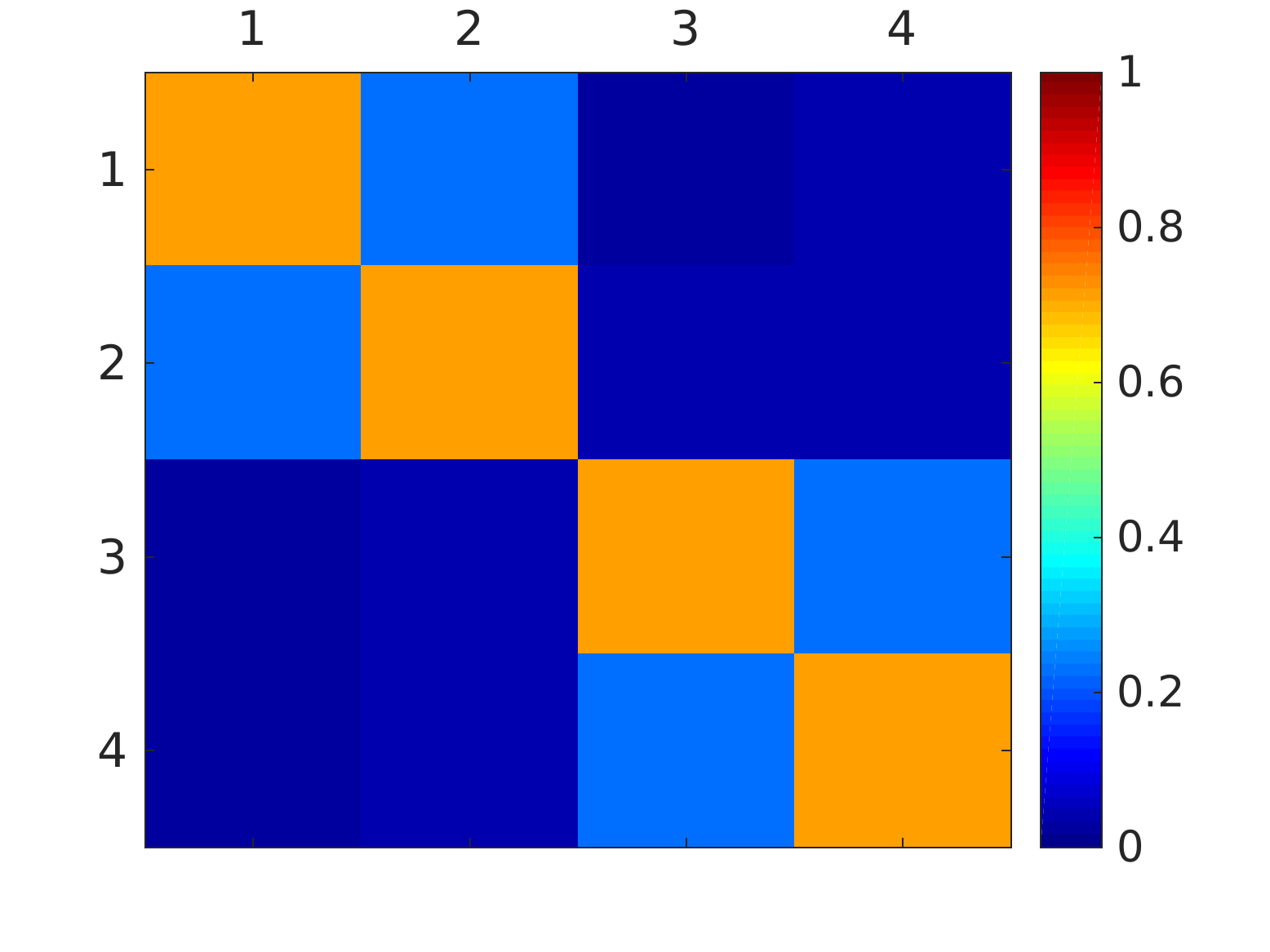}
		(b) Strongly metastable matrix $P_c = S^{-1} T$.
	\end{minipage}
	\captionof{figure}{Coupling matrix and projected transition matrix for a lag-time $\tau_1 = 0.2$.}
	\label{fig:electron_stability1}
\end{minipage} \\

\begin{minipage}{\linewidth}
	\centering
	\begin{minipage}[t]{0.29\textwidth}
		\includegraphics[width=\textwidth]{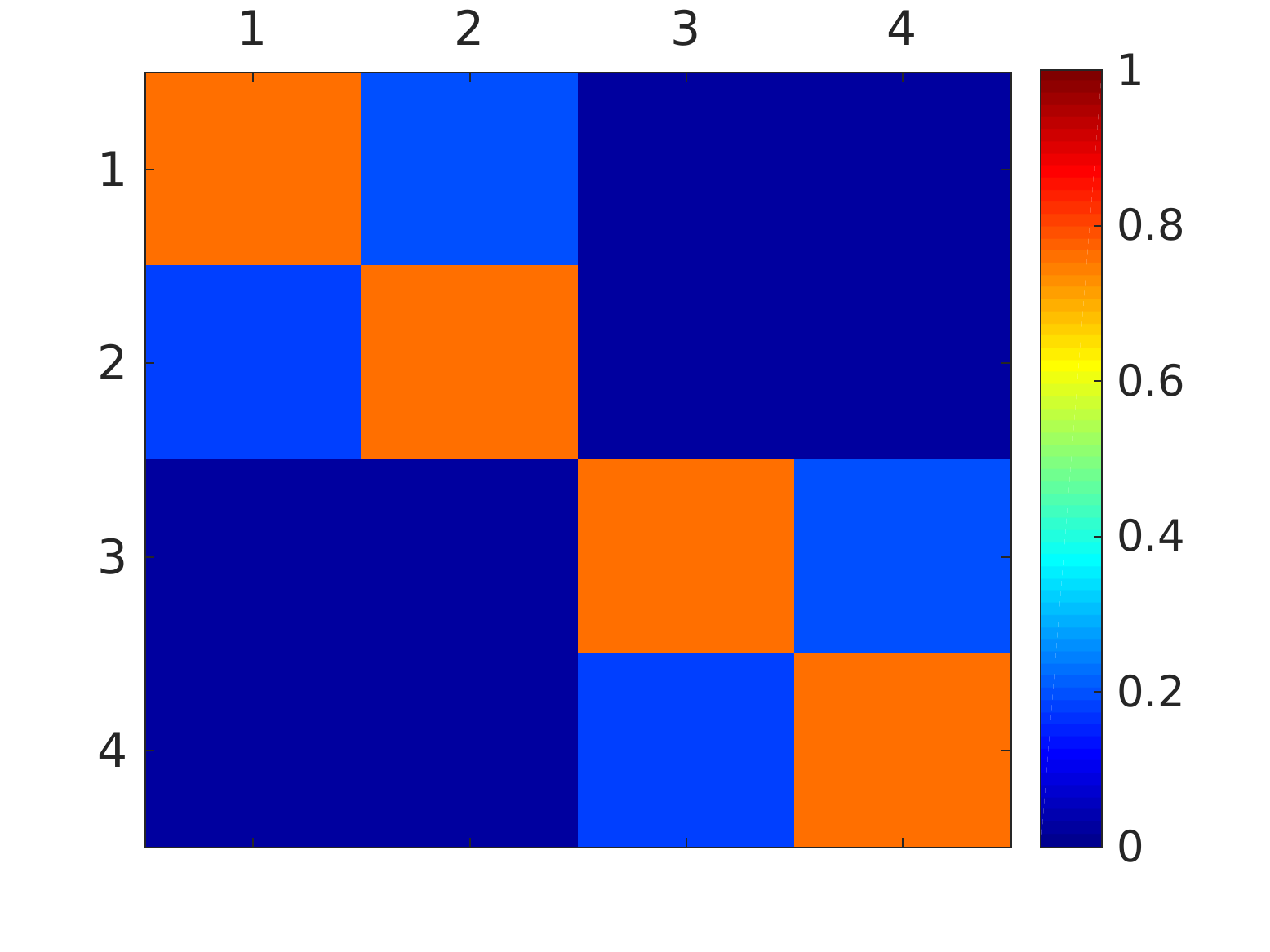}
		(a) Weakly metastable coupling matrix $T$.
	\end{minipage}
	\hspace{20pt}
	\begin{minipage}[t]{0.29\textwidth}
		\includegraphics[width=\textwidth]{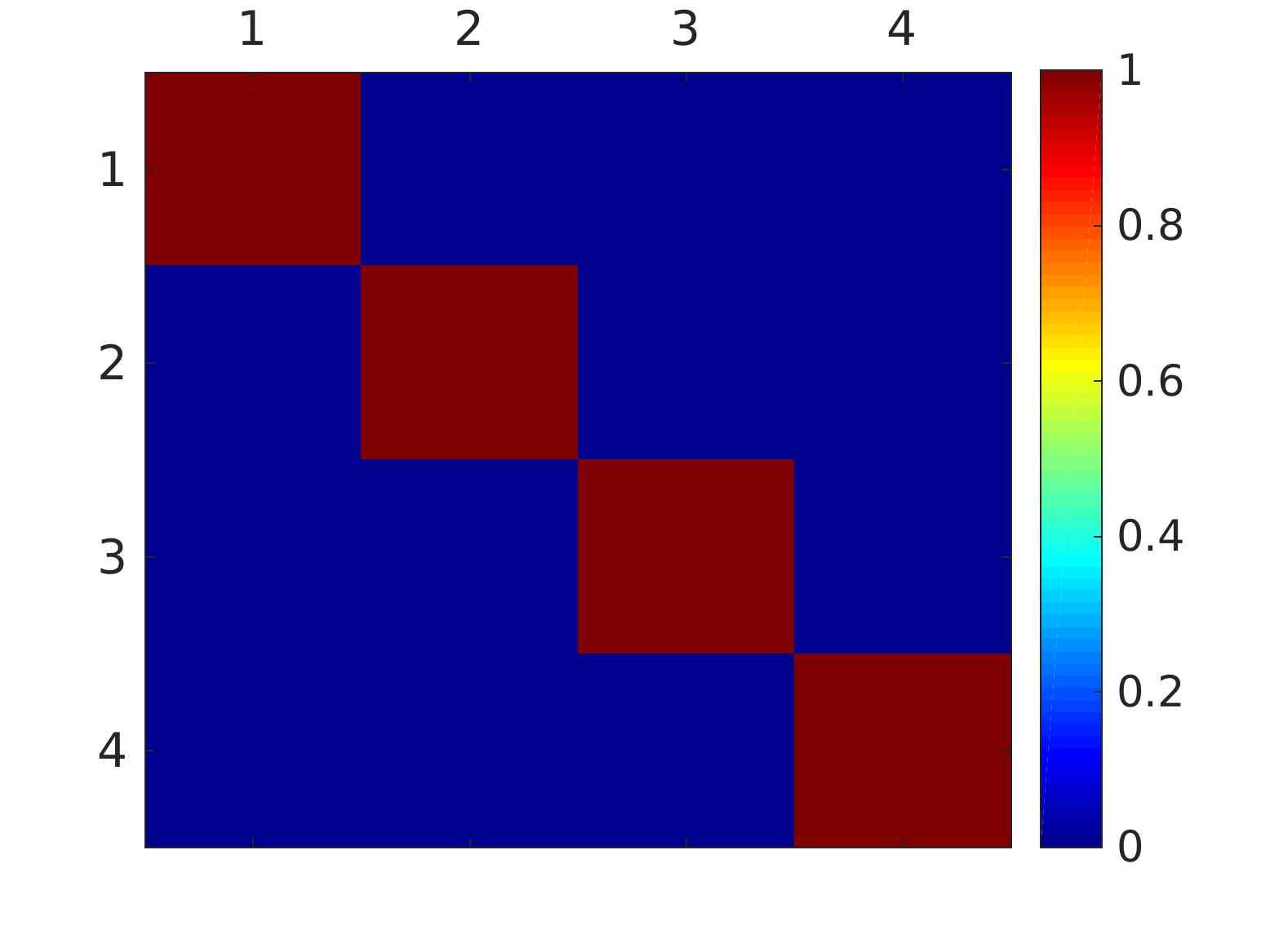}
		(b) Strongly metastable matrix $P_c = S^{-1} T$.
	\end{minipage}
	\captionof{figure}{Coupling matrix and projected transition matrix for a small lag-time $\tau_2 = 10^{-3}$.}
	\label{fig:electron_stability2}
\end{minipage} \\

\begin{table}[t!]
	\centering 
	\begin{tabular}{|c|c|c|c|c|} \hline
		Metastable subset & $1$ & $2$  & $3$ & $4$ \\ \hline \hline
		Statistical weight & $0.2406$ & $0.2556$ & $0.2520$ & $0.2518$ \\ \hline \hline
		Metastability $T (\tau_1)$ & $0.5811$  & $0.5827$ & $0.5884$ & $0.5815$ \\ \hline
		Metastability $P_c (\tau_1)$ & $0.7077$  & $0.7084$ & $0.7135$ & $0.7082$  \\ \hline \hline
		Metastability $T (\tau_2)$ & $0.7571$  & $0.7577$ & $0.7622$ & $0.7577$ \\ \hline
		Metastability $P_c (\tau_2)$ & $0.9980$  & $0.9980$ & $0.9980$ & $0.9980$  \\ \hline
	\end{tabular}
	\caption{Influence of rebinding to the stability of $P_c$ for different lag-times $\tau_1, \tau_2$.}
\end{table}

\subsection{An almost reversible process}\label{sec:addit}
Only the matrix $Q_c$ is needed to estimate the minimal rebinding effect. If this matrix is reversible, then the minimal rebinding effect cannot be effectively provided, as the previous section \ref{sec:ED} has shown. Thus, non-reversibility of $Q_c$ is the crucial prerequisite for the estimation of the rebinding effect which has also been demonstrated in Sec. \ref{sec:art}. There are examples where the detailed process is (almost) reversible, but the projection $Q_c$ is not. Such an example will be discussed now.  

\begin{figure}[ht!]
	\centering
	\includegraphics[width=0.4\textwidth]{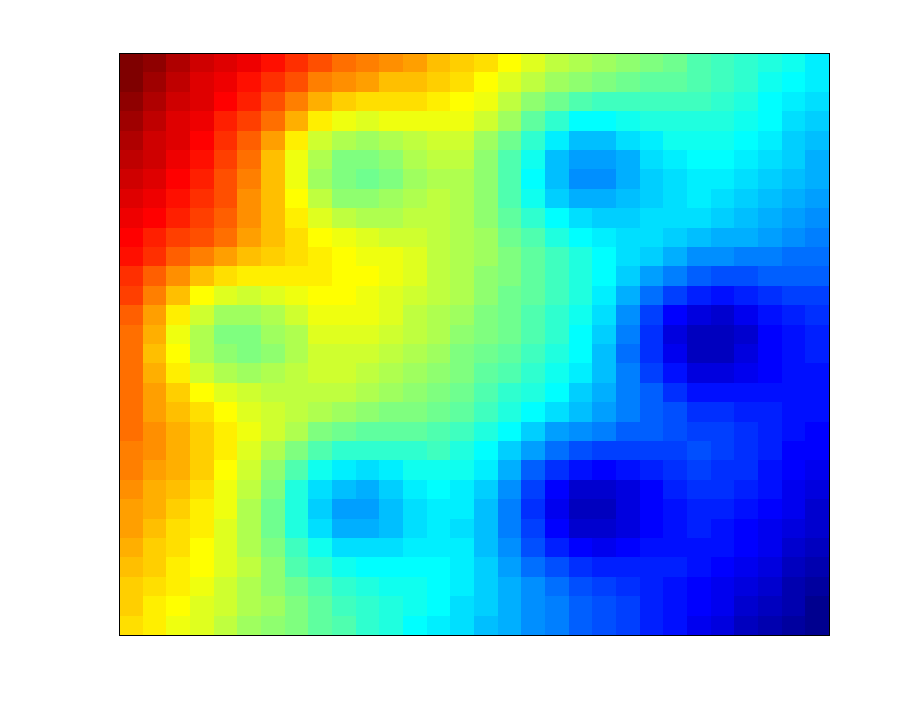}
	\caption{From the invariant density of the manipulated SQRA-process, one can compute a free energy landscape (taking the negative logarithm of the entries of the invariant density vector). One can clearly see, that this energy landscape is tilted to the right lower corner of the state space.}
	\label{fig:freeEnergy}
\end{figure}

We consider an invariant density with 6 Gaussians which are arranged in a circle. The domain is
then decomposed into $30 \times 30$ boxes (grid).
In each center point $p_i$ of the 900 Boxes $\{B_i\}_{i=1,...,900}$
we evauated the corresponding invariant density $\pi_i$.
The transition rates $q_{ij}$ of the transition matrix $Q\in \mathbb{R}^{900\times 900}$ were computed according to the square root approximation \cite{SquareRootApprox}, i.e.
\[
q_{ij}=\sqrt{\frac{\pi_j}{\pi_i}}
\]
for neighbouring boxes $B_i$ and $B_j$.
We remark that there also exist other methods to assemble the rate matrix $Q$, for instance the milestoning method \cite{milestoning}, where not the full state space but only the part with the metastabilities is discretized.
In order to provide a non-reversible process, all transitions which go from a box with a smaller index  to a box with a higher index are multiplied with a factor $1.2$.
Using the GenPCCA method the $900\times 900$ matrix is projected onto a $3\times 3$-matrix $Q_c$. 
The resulting projection is:
\begin{equation*}
Q_c = 
\begin{pmatrix}
   -0.0263 &   0.0219 &   0.0044\\
    0.0025 &  -0.0174 &   0.0149\\
    0.0022 &   0.0195 &  -0.0217
\end{pmatrix}.
\end{equation*}

Based on the results of GenPPCA, with membership functions $\chi=XA$, the real rebinding effect is given by $\mathrm{det}(S_{real})=0.0031$. In order to estimate this effect from the matrix $Q_c$, we first compute the eigenvalues of $Q_c$ which are real valued. Thus, the pattern of the Schur matrix $\Xi$ depends on our decision whether the original process is reversible or not. If we assume, that the original process is non-reversible (which is indeed the case), then only $\clubsuit$ in (\ref{eq:schur_structure}) is zero, otherwise, $\clubsuit$ and $\spadesuit$ are supposed to be zero.  

In order to maximize the determinant of the matrix $S$ with constraints $A\in{\cal C}$, the global optimization problem has been solved by a multi-start ansatz (500 starts) and with a quadratic penalty function approach.  

In the non-reversible setting we find an optimal linear transformation matrix $A$ leading to a rebinding effect of $\mathrm{det}(S_{opt})=0.7446$. The original process $Q$ is almost reversible. Assuming a reversible original process further restricts the set $\cal C$ of feasible transformation matrices, the estimate becomes ``better'' (lower determinant). In this case, the optimal linear transformation matrix $A$ leads to an estimated rebinding effect of $\mathrm{det}(S_{opt})=0.3179$. In both cases, the correct pattern of $\Xi$ with the eigenvalues of $Q_c$ on its diagonal has been revealed. In the non-reversible case we get $\spadesuit=0.0036$ (the true value is $0.0001$, we overestimated the non-reversibility of the original process). In the reversible case we get by construction $\spadesuit=0$. Thus, the knowledge about the reversibility or non-reversibility of the original process can improve the estimate of the rebinding effect a lot. In other words, computing the minimal possible rebinding effect of non-reversible processes provides just a rough estimate as Sec. \ref{sec:ED} and Sec. \ref{sec:addit} have demonstrated. 

\begin{figure}[ht!]
	\centering
	\includegraphics[width=0.42\textwidth]{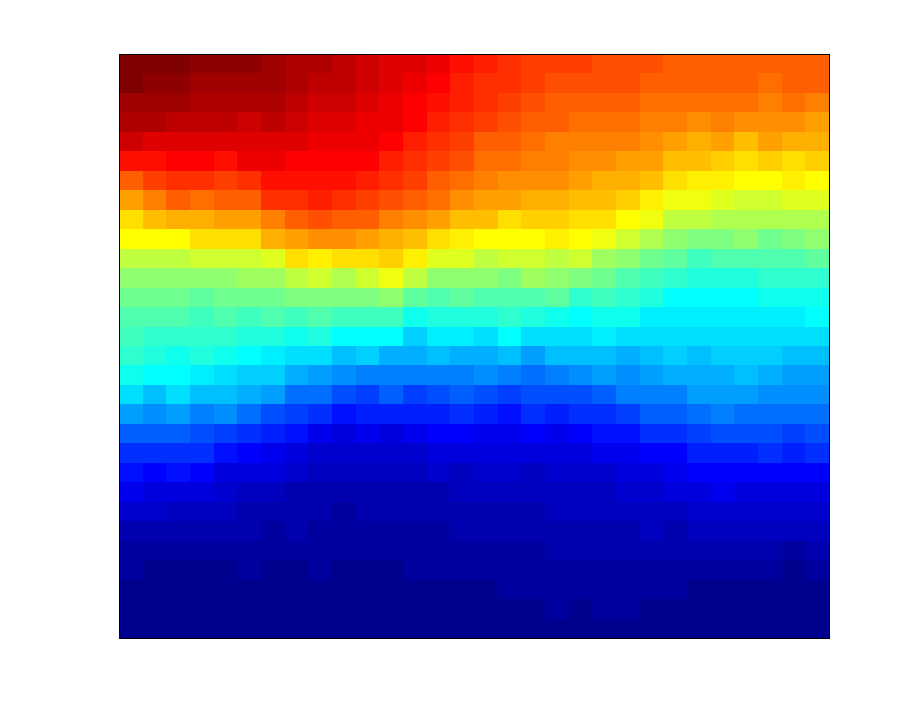}
	\includegraphics[width=0.4\textwidth]{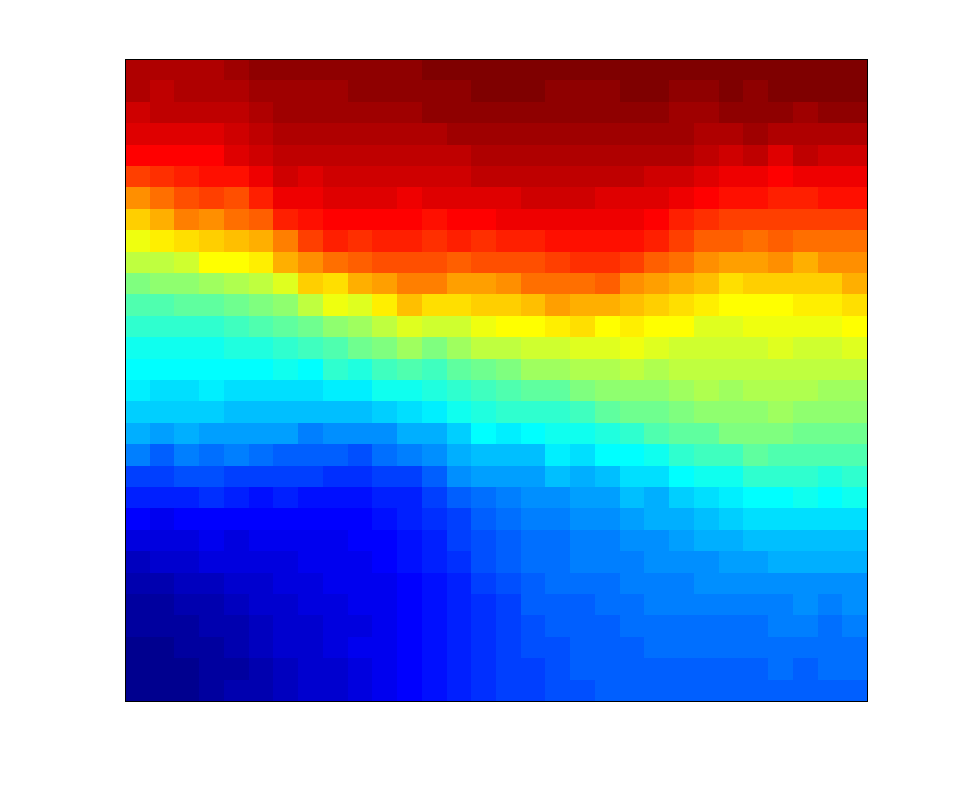}
	\caption{Left: One of the three membership functions obtained by GenPCCA applied to the original process. Right: Applying the optimized transformation matrix (the result of the optimization problem) to the original Schur vectors also leads to a partition of unity set of membership functions $\chi_{opt}$. The shape of these functions are similar to those of the original process, but they are not feasible membership functions anymore. The plotted function has values between $-0.8735$ and $1.4020$ and has a similar shape compared to the left one.}
	\label{fig:membership_compare}
\end{figure}

\section{Conclusion}
\label{sec:conclusion}
In this paper, two recent research topics were combined by extending the computation of a lower bound for the rebinding effect onto non-reversible processes.
The generalized fuzzy clustering algorithm GenPCCA has been employed to obtain the optimal membership functions as a linear combination of the dominant Schur vectors.
The overlap of the membership functions is crucial for a correct mapping, though influences the observed stability of the system. The more overlap, the more stable the macro states appear to be. \\
This phenomenon is denoted as \textbf{rebinding effect} because of its occurrence in receptor-ligand-systems, where this `spatial memory' leads to an increased probability for a fast rebinding after the dissociation of a receptor-ligand-complex. Under the assumption of a fuzzy clustering $\chi = XA$, the minimal rebinding effect included in a given kinetics has been computed as the solution of an optimization problem, considering reversible as well as non-reversible processes by using Schur vectors $X$. This optimization problem has been tested for some numerical examples, showing that the quality of the estimation can be good or bad, yet becomes less predictable for large degrees of non-reversibility of $Q_c$. \\

Knowing the rebinding effect of a system can be of particular relevance for applications like computational drug design, where it is essential to correctly predict binding affinities in order to evaluate the expected efficiency of a newly designed drug. Since many real-world processes are non-reversible, it was important to add this case to the already existing optimization problem for reversible processes. This extension yields an estimation for the rebinding effect of a clustered system, without the necessity to know if the original process was actually reversible or non-reversible. \\

In this paper, the rebinding effect has been tackled from a rather theoretical perspective. For further research, it could be of interest to combine and extend the obtained results with the outcomes from molecular dynamics simulations.

A first approach towards this method can be found in \cite{Roehl2017}.
However this article now, contains substantial changes.



\section*{Acknowledgments}
This work was supported by Math+ and by the 
CRC-1114 'Scaling Cascades in Complex Systems', project A05.







\end{document}